\documentclass[a4paper, onecolumn,10pt]{article}
\usepackage[english]{babel}
\usepackage[T1]{fontenc}
\usepackage[latin1]{inputenc}
\usepackage{graphicx}

\usepackage{pgfplots}
\usepackage{pgfplotstable}
\usepackage{colortbl}

\usepgfplotslibrary{external}
\tikzsetexternalprefix{pgfplotsfigures/}
\tikzset{external/optimize=false}

\usepackage{graphicx}

\usepackage[small,bf]{caption}

\usepackage{enumitem}

\usepackage[normal,tight]{subfigure}
\newcommand{\goodgap}{%
\hspace{\subfigtopskip}%
\hspace{\subfigbottomskip}}
\makeatletter
  \def\zaptype#1{%
  \listsubcaptions 
  \def\@captype{#1}}
\makeatother

\usepackage{amsmath}
\usepackage{amssymb}
\usepackage{amsfonts}
\usepackage{amsthm}
\theoremstyle{definition}

\newtheorem{corollary}{Corollary}
\newtheorem{definition}{Definition}
\newtheorem{lemma}{Lemma}
\newtheorem{proposition}{Proposition}
\newtheorem{theorem}{Theorem}

\newtheorem{example}{Example}
\theoremstyle{remark}

\newenvironment{keywords}{\begin{quote}%
  \textbf{Keywords: }}{\hfill
\end{quote}%
}
\newenvironment{AMS}{\begin{quote}%
  \textbf{AMS: }}{\hfill
\end{quote}%
}

\usepackage{ifthen}

\makeatletter
\def\s@btitle{\relax}
\def\subtitle#1{\gdef\s@btitle{#1}}
\def\@maketitle{%
  \newpage
  \null
  \vskip 2em%
  \begin{center}%
  \let \footnote \thanks
    {\LARGE \@title \par}%
                \if\s@btitle\relax
                \else\typeout{[subtitle]}%
                        \vskip .5pc
                        \begin{large}%
                                \textsl{\s@btitle}%
                                \par
                        \end{large}%
                \fi
    \vskip 1.5em%
    {\large
      \lineskip .5em%
      \begin{tabular}[t]{c}%
        \@author
      \end{tabular}\par}%
    \vskip 1em%
    {\large \@date}%
  \end{center}%
  \par
  \vskip 1.5em}
\makeatother 

\usepackage{color}
\definecolor{DarkGray}{rgb}{0.1,0.1,0.1}
\definecolor{BlueViolet}{rgb}{0.2,0.,0.4}
\definecolor{Grey}{rgb}{0.3,0.3,0.3}
\definecolor{Red}{rgb}{0.8,0.,0.}
\definecolor{Yellow}{rgb}{0.,0.4,0.4}
\definecolor{RoyalBlue}{rgb}{0.,0.14,0.8}
\definecolor{BrickRed}{rgb}{0.8,0.34,0.34}
\definecolor{LimeGreen}{rgb}{0.,0.65,0.25}
\definecolor{Sepia}{rgb}{0.4,0.3,0.3}

\usepackage{hyperref}
\hypersetup{
    colorlinks,
    citecolor=BrickRed,
    filecolor=LimeGreen,
    linkcolor=RoyalBlue,
    urlcolor=Sepia
}

\usepackage{accents}

\usepackage{dsfont}

\def\reel{\mathbb{R}}
\def\nat{\mathbb{N}}

\def\citechapter{Chap.}
\def\citesection{Sect.}
\def\citepage{p.}

\DeclareMathOperator{\Span}{span}
\DeclareMathOperator{\argmin}{argmin}
\providecommand{\abs}[1]{\left\lvert#1\right\rvert}
\providecommand{\norm}[1]{\left\lVert#1\right\rVert}
\def\adj{\ast}

\newcommand{\dual}[1]{\ensuremath{{#1}^{\adj}}}%
\newcommand{\DualSpace}[1]{\ensuremath{{\cal #1}^{\adj}}}%
\newcommand{\skp}[2]{\ensuremath{\langle #1 , #2 \rangle}}%
\newcommand{\SKP}[2]{\ensuremath{\left \langle #1 , #2 \right \rangle}}%
\newcommand{\Vect}[1]{\ensuremath{{\bf #1}}}
\newcommand{\Mat}[1]{\ensuremath{{\bf #1}}}
\newcommand{\Op}[1]{\ensuremath{\hat{#1}}}

\newcommand{\Space}[1]{{\cal #1}}
\newcommand{\Null}[1]{\textrm{nul}\left ( #1 \right ) }
\newcommand{\Range}[1]{\textrm{ran}\left ( #1 \right ) }

\usepackage[acronym,translate=false,counter=equation,sanitize={name=false,description=true,symbol=false}]{glossaries}
\glsdisablehyper

\newglossary[dlg]{definitionlist}{dyg}{def}{Definitions}
\newglossary[slg]{symbolslist}{syg}{sym}{Symbols}

\usepackage{ifnextok}

\def\sym{\sympage}

\def\glosentryexists#1{%
  \ifglsentryexists{def:#1}{%
    \def\myglossaryprefix{def}%
    \ifglsentryexists{glos:#1}{%
      \glsadd[format=hyperbf,counter=page]{glos:#1}%
    }{%
      \relax%
    }%
  }{%
    \def\myglossaryprefix{glos}%
  }%
}

\makeatletter
\def\glosname#1{%
\glosentryexists{#1}%
\IfNextToken\bgroup%
{\glosnametwoargs@aux{#1}}{\glosnameonearg@only{#1}}%
}%
\def\glosnameonearg@only#1{\glsname[format=hyperbf,counter=page]{\myglossaryprefix:#1}}%
\def\glosnametwoargs@aux#1#2{\glsname[format=hyperbf,counter=page]{\myglossaryprefix:#1}[#2]}%
\makeatother

\makeatletter
\def\glostext#1{%
\glosentryexists{#1}%
\IfNextToken\bgroup%
{\glostexttwoargs@aux{#1}}{\glostextonearg@only{#1}}%
}%
\def\glostextonearg@only#1{\gls[format=hyperbf,counter=page]{\myglossaryprefix:#1}}%
\def\glostexttwoargs@aux#1#2{\gls[format=hyperbf,counter=page]{\myglossaryprefix:#1}[#2]}%
\makeatother

\makeatletter
\def\glospl#1{%
\glosentryexists{#1}%
\IfNextToken\bgroup%
{\glospltwoargs@aux{#1}}{\glosplonearg@only{#1}}%
}%
\def\glosplonearg@only#1{\glspl[format=hyperbf,counter=page]{\myglossaryprefix:#1}}%
\def\glospltwoargs@aux#1#2{\glspl[format=hyperbf,counter=page]{\myglossaryprefix:#1}[#2]}%
\makeatother

\makeatletter
\def\glosfirst#1{%
\glosentryexists{#1}%
\IfNextToken\bgroup%
{\glosfirsttwoargs@aux{#1}}{\glosfirstonearg@only{#1}}%
}%
\def\glosfirstonearg@only#1{\glsfirst[format=hyperbf,counter=page]{\myglossaryprefix:#1}}%
\def\glosfirsttwoargs@aux#1#2{\glsfirst[format=hyperbf,counter=page]{\myglossaryprefix:#1}[#2]}%
\makeatother

\makeatletter
\def\glosfirstpl#1{%
\glosentryexists{#1}%
\IfNextToken\bgroup%
{\glosfirstpltwoargs@aux{#1}}{\glosfirstplonearg@only{#1}}%
}%
\def\glosfirstplonearg@only#1{\glsfirstplural[format=hyperbf,counter=page]{\myglossaryprefix:#1}}%
\def\glosfirstpltwoargs@aux#1#2{\glsfirstplural[format=hyperbf,counter=page]{\myglossaryprefix:#1}[#2]}%
\makeatother

\makeatletter
\def\Glosname#1{%
\glosentryexists{#1}%
\IfNextToken\bgroup%
{\Glosnametwoargs@aux{#1}}{\Glosnameonearg@only{#1}}%
}%
\def\Glosnameonearg@only#1{\Gls[format=hyperbf,counter=page]{\myglossaryprefix:#1}}%
\def\Glosnametwoargs@aux#1#2{\Gls[format=hyperbf,counter=page]{\myglossaryprefix:#1}[#2]}%
\makeatother

\makeatletter
\def\Glospl#1{%
\glosentryexists{#1}%
\IfNextToken\bgroup%
{\Glospltwoargs@aux{#1}}{\Glosplonearg@only{#1}}%
}%
\def\Glosplonearg@only#1{\Glspl[format=hyperbf,counter=page]{\myglossaryprefix:#1}}%
\def\Glospltwoargs@aux#1#2{\Glspl[format=hyperbf,counter=page]{\myglossaryprefix:#1}[#2]}
\makeatother

\makeatletter
\def\GLOSname#1{%
\glosentryexists{#1}%
\IfNextToken\bgroup%
{\GLOSnametwoargs@aux{#1}}{\GLOSnameonearg@only{#1}}%
}%
\def\GLOSnameonearg@only#1{\GLS[format=hyperbf,counter=page]{\myglossaryprefix:#1}}%
\def\GLOSnametwoargs@aux#1#2{\GLS[format=hyperbf,counter=page]{\myglossaryprefix:#1}[#2]}%
\makeatother

\makeatletter
\def\GLOSpl#1{%
\glosentryexists{#1}%
\IfNextToken\bgroup%
{\GLOSpltwoargs@aux{#1}}{\GLOSplonearg@only{#1}}%
}%
\def\GLOSplonearg@only#1{\GLSpl[format=hyperbf,counter=page]{\myglossaryprefix:#1}}%
\def\GLOSpltwoargs@aux#1#2{\GLSpl[format=hyperbf,counter=page]{\myglossaryprefix:#1}[#2]}%
\makeatother


\def\dummysympage{\sympage}
\makeatletter
\def\glosformula#1{%
  \ifx\sym\dummysympage%
    \glosformula@nonnumbered{#1}%
  \else%
    \glosformula@numbered{#1}%
  \fi%
}%
\def\glosformula@numbered#1{%
\glsentryuseri{\myglossaryprefix:#1}%
}%
\def\glosformula@nonnumbered#1{%
\glsentryuseri{\myglossaryprefix:#1}%
}%
\makeatother

\newenvironment{myequation*}{\def\sym{\sympage}\equation\aligned\protect}{\nonumber\endaligned\endequation}%
\makeatletter
\newcommand\mynonumber{\start@align\@ne\st@rredtrue\m@ne}%
\newcommand\mynumber{\start@align\@ne\st@rredfalse\m@ne}%
\makeatother

\newenvironment{myalign*}{\def\sym{\sympage}\mynonumber\protect}{\endalign}%
\newcommand{\symeq}[1]{
  \protect\gls[format=hyperit]{sym:#1}
  \ifglsentryexists{def:#1}{%
    \glsadd[format=hyperbf,counter=page]{def:#1}
  }{%
    \relax%
  }
}
\newcommand{\sympage}[1]{\glssymbol[format=hyperbf,counter=page]{def:#1}}

\newcommand{\linkonly}[1]{\glsadd[counter=page]{#1}\glshyperlink{#1}}

\newglossarystyle{definitiontable}{%
{\begin{longtable}{cp{5cm}p{6cm}p{\glspagelistwidth}}}%
{\end{longtable}}%
\renewcommand*{\glsgroupheading}[1]{}%
\renewcommand*{\glossaryentryfield}[5]{%
##4
& \itshape{\glstarget{##1}{##2}}
&%
\ifthenelse{\equal{\glsentryuseri{##1}}{\empty}}
{%
\relax%
}{%
\begin{equation}\label{math-##1}\glsentryuseri{##1}\end{equation}
}%
& ##5
\\
& 
\multicolumn{2}{p{11cm}}{##3}
\\
}%
\renewcommand*{\glossarysubentryfield}[6]{%
##5
& \itshape{\glstarget{##2}{##3}}
& 
\ifthenelse{\equal{\glsentryuseri{##2}}{\empty}}%
{%
\relax%
}{%
\begin{equation}\label{math-##2}\glsentryuseri{##2}\end{equation}
}%
& ##6
\\
& 
\multicolumn{2}{p{11cm}}{##4}
\\
}
}

\makeglossaries

\newtheorem{method}{Method}

\newglossaryentry{glos:ALS}{
	type=\acronymtype,
	name={Alternating Least Squares},
	text={ALS},
	first={Alternating Least Squares (ALS)},
	description={In matrix factorization alternatingly one of the matrix factors is kept fixed and the residual minimized in the $\ell_2$-norm with respect to the other factor},
	sort=ALS
}
\newglossaryentry{glos:BSS}{
	type=\acronymtype,
	name={Blind Source Separation},
	text={BSS},
	first={Blind Source Separation (BSS)},
	description={In the problem of Blind Source Separation a given generating system is used to best represent a given signal under given constraints such as sparsity of the representation},
	sort=BSS
}
\newglossaryentry{glos:CG}{
	type=\acronymtype,
	name={Conjugate Gradient},
	text={CG},
	first={Conjugate Gradient (CG)},
	description={Using the conjugacy of consecutive search directions possible in Hilbert spaces, the CG method shows superlinear behavior, overcoming the Steepest Descent methods problems with ill-conditioned matrices},
	sort=CG
}
\newglossaryentry{glos:CGNE}{
	type=\acronymtype,
	name={Conjugate Gradient Normal Error},
	text={CGNE},
	first={Conjugate Gradient Normal Error (CGNE)},
	description={Similar to the \glostext{CG} method, this is a modification for the normal equation to end up with an the error normal to the residual},
	sort=CGNE
}
\newglossaryentry{glos:CGNR}{
	type=\acronymtype,
	name={Conjugate Gradient Normal Residual},
	text={CGNR},
	first={Conjugate Gradient Normal Residual (CGNR)},
	description={Similar to the \glostext{CG} method, this is a modification for the normal equation with a normal residual},
	sort=CGNR
}
\newglossaryentry{glos:RESESOP}{
	type=\acronymtype,
	name={Regularized Sequential Subspace Optimization},
	text={RESESOP},
	first={Regularized Sequential Subspace Optimization (RESESOP)},
	description={Generalization of the Landweber algorithm to multiple search directions using Bregman projections for noisy data},
	sort=RESESOP
}
\newglossaryentry{glos:SD}{
	type=\acronymtype,
	name={Steepest Descent},
	text={SD},
	first={Steepest Descent (SD)},
	description={SD uses the gradient as descent direction and performs a one-dimensional line search along it to converge to the minimum},
	sort=SD
}
\newglossaryentry{glos:SVD}{
	type=\acronymtype,
	name={Singular Value Decomposition},
	text={SVD},
	first={Singular Value Decomposition (SVD)},
	description={For each matrix a unique decomposition into the singular triplets $\{\sigma_i, u_i, v_i\}$ exists such that the original matrix is the product $U \Sigma V$, where $\Sigma$ is the matrix containing the singular values $\sigma_i$ on the diagonal. It is similar to the eigendecomposition for symmetric matrices},
	sort=SD
}
\newglossaryentry{glos:SESOP}{
	type=\acronymtype,
	name={Sequential Subspace Optimization},
	text={SESOP},
	first={Sequential Subspace Optimization (SESOP)},
	description={Generalization of the Landweber algorithm to multiple search directions using Bregman projections},
	sort=SESOP
}

 \newglossaryentry{def:dimension}{
	type=definitionlist,
	symbol={\ensuremath{n}},
 	name={discretization dimension},
 	description={The number of points in the interval $[0,1]$ used to discretize the continuous problem},
 	sort=A
 }
 \newglossaryentry{def:distance-bregman}{
	type=definitionlist,
	symbol={\ensuremath{\Delta}},
 	name={Bregman distance},
 	description={The Bregman distance is a generalization of the euclidian distances for Banach spaces},
 	sort=Delta
 }
 \newglossaryentry{def:mapping-duality}{
	type=definitionlist,
	symbol={\ensuremath{J}},
 	name={duality mapping},
 	description={},
 	sort=J
 }
 \newglossaryentry{def:mapping-duality-single-valued}{
	type=definitionlist,
	symbol={\ensuremath{j}},
 	name={single-valued duality mapping},
 	description={obtained by a specific selection from the set defined by the duality mapping},
 	sort=j
 }
 \newglossaryentry{def:matrix-forward}{
	type=definitionlist,
	symbol={\ensuremath{\Mat{A}}},
 	name={forward matrix},
 	description={The forward matrix is the central element of the discretized forward problem $\Mat{A}x = y$, where we have a vector $x$ and we look for a vector $y$},
 	sort=A
 }
 \newglossaryentry{def:operator-forward}{
	type=definitionlist,
	symbol={\ensuremath{\protect\Op{A}}},
 	name={forward operator},
 	description={The forward operator is the central element of the forward problem $\Mat{A}x = y$, where we have $x$ and we look for $y$},
 	sort=A
 }
 \newglossaryentry{def:projection-bregman}{
	type=definitionlist,
	symbol={\ensuremath{\Pi}},
 	name={Bregman projection},
 	description={The Bregman projection is a generalization of the metric projection for Banach spaces},
 	sort=Pi
 }
 \newglossaryentry{def:projection-metric}{
	type=definitionlist,
	symbol={\ensuremath{P}},
 	name={Metric projection},
 	description={The metric projection with respect to a convex set is one point in the set that minimizes the norm distance to the given point},
 	sort=P
 }

\newglossaryentry{def:dimension-rhs}{
	type=definitionlist,
	symbol={\ensuremath{l}},
 	name={discretization dimension for right-hand side},
 	description={The number of points in the interval $[0,1]$ used to discretize the right-hand side function},
 	sort=l
}
\newglossaryentry{def:dimension-solution}{
	type=definitionlist,
	symbol={\ensuremath{m}},
 	name={discretization dimension for solution},
 	description={The number of points in the interval $[0,1]$ used to discretize the solution function},
 	sort=m
}
\newglossaryentry{def:iterations}{
	type=definitionlist,
	symbol={\ensuremath{I}},
 	name={number of iterations},
 	description={Number of iterations to perform the minimization},
 	sort=I
}
\newglossaryentry{def:iterations-limit}{
	type=definitionlist,
	symbol={\ensuremath{I^{\text{max}}}},
 	name={maximum number of iterations},
 	description={Threshold for the number of iterations to perform the minimization},
 	sort=Imax
}
\newglossaryentry{def:norm-dualX}{
	type=definitionlist,
	symbol={\ensuremath{\dual{p}}},
 	name={Lp norm of $\DualSpace{X}$},
 	description={Value $\dual{p}$ of the Lp norm of $\DualSpace{X}$},
 	sort=p
}
\newglossaryentry{def:norm-X}{
	type=definitionlist,
	symbol={\ensuremath{p}},
 	name={Lp norm of $\Space{X}$},
 	description={Value $p$ of the Lp norm of $\Space{X}$, dual to the norm of $\Space{X}$},
 	sort=q
}
\newglossaryentry{def:norm-Y}{
	type=definitionlist,
	symbol={\ensuremath{r}},
 	name={Lp norm of $\Space{Y}$},
 	description={Value $p$ of the Lp norm of $\Space{Y}$},
 	sort=r
}
\newglossaryentry{def:number-search-directions}{
	type=definitionlist,
	symbol={\ensuremath{N}},
 	name={number of search directions},
 	description={The number of search directions is basically the dimension of the search space $\Space{U}$},
 	sort=N
}
\newglossaryentry{def:power-X}{
	type=definitionlist,
	symbol={\ensuremath{\pi_{\Space{X}}}},
 	name={power type in $\Space{X}$},
 	description={power type of a duality mapping from $\Space{X}$ to $\DualSpace{X}$},
 	sort=powerx
}
\newglossaryentry{def:power-Y}{
	type=definitionlist,
	symbol={\ensuremath{\pi_{\Space{Y}}}},
 	name={power type in $\Space{Y}$},
 	description={power type of a duality mapping from $\Dual{\Space{Y}}$ to $\Space{Y}$},
 	sort=powery
}
\newglossaryentry{def:runtime}{
	type=definitionlist,
	symbol={\ensuremath{T}},
 	name={runtime},
 	description={Required walltime for the minimization function},
 	sort=T
}
\newglossaryentry{def:threshold-residual}{
	type=definitionlist,
	symbol={\ensuremath{\Delta}},
 	name={residual threshold},
 	description={Iteration is stopped if residual drops below this threshold in magnitude},
 	sort=delta
}

 \newglossaryentry{sym:dimension}{
 	type=symbolslist,
 	name={\glsentrysymbol{def:dimension}},
 	description={\linkonly{def:dimension}},
	sort={A}
 }
 \newglossaryentry{sym:distance-bregman}{
 	type=symbolslist,
 	name={\glsentrysymbol{def:distance-bregman}},
 	description={\linkonly{def:distance-bregman}},
	sort={Delta}
 }
 \newglossaryentry{sym:mapping-duality}{
 	type=symbolslist,
 	name={\glsentrysymbol{def:mapping-duality}},
 	description={\linkonly{def:mapping-duality}},
	sort={J}
 }
 \newglossaryentry{sym:mapping-duality-single-valued}{
 	type=symbolslist,
 	name={\glsentrysymbol{def:mapping-duality-single-valued}},
 	description={\linkonly{def:mapping-duality-single-valued}},
	sort={j}
 }
 \newglossaryentry{sym:matrix-forward}{
 	type=symbolslist,
 	name={\glsentrysymbol{def:matrix-forward}},
 	description={\linkonly{def:matrix-forward}},
	sort={A}
 }
 \newglossaryentry{sym:operator-forward}{
 	type=symbolslist,
 	name={\glsentrysymbol{def:operator-forward}},
 	description={\linkonly{def:operator-forward}},
	sort={A}
 }
 \newglossaryentry{sym:projection-bregman}{
 	type=symbolslist,
 	name={\glsentrysymbol{def:projection-bregman}},
 	description={\linkonly{def:projection-bregman}},
	sort={Delta}
 }

\newglossaryentry{sym:dimension-rhs}{
	type=symbolslist,
 	name={\glsentrysymbol{def:dimension-rhs}},
 	description={\linkonly{def:dimension-rhs}},
	sort={l}
}
\newglossaryentry{sym:dimension-solution}{
	type=symbolslist,
 	name={\glsentrysymbol{def:dimension-solution}},
 	description={\linkonly{def:dimension-solution}},
	sort={m}
}
\newglossaryentry{sym:iterations}{
	type=symbolslist,
 	name={\glsentrysymbol{def:iterations}},
 	description={\linkonly{def:iterations}},
	sort={I}
}
\newglossaryentry{sym:iterations-limit}{
	type=symbolslist,
 	name={\glsentrysymbol{def:iterations-limit}},
 	description={\linkonly{def:iterations-limit}},
	sort={Imax}
}
\newglossaryentry{sym:norm-dualX}{
	type=symbolslist,
 	name={\glsentrysymbol{def:norm-dualX}},
 	description={\linkonly{def:norm-dualX}},
	sort={q}
}
\newglossaryentry{sym:norm-X}{
	type=symbolslist,
 	name={\glsentrysymbol{def:norm-X}},
 	description={\linkonly{def:norm-X}},
	sort={p}
}
\newglossaryentry{sym:norm-Y}{
	type=symbolslist,
 	name={\glsentrysymbol{def:norm-Y}},
 	description={\linkonly{def:norm-Y}},
	sort={r}
}
\newglossaryentry{sym:number-search-directions}{
	type=symbolslist,
 	name={\glsentrysymbol{def:number-search-directions}},
 	description={\linkonly{def:number-search-directions}},
	sort={N}
}
\newglossaryentry{sym:power-X}{
	type=symbolslist,
 	name={\glsentrysymbol{def:power-X}},
 	description={\linkonly{def:power-X}},
	sort={pix}
}
\newglossaryentry{sym:power-Y}{
	type=symbolslist,
 	name={\glsentrysymbol{def:power-Y}},
 	description={\linkonly{def:power-Y}},
	sort={piy}
}
\newglossaryentry{sym:runtime}{
	type=symbolslist,
 	name={\glsentrysymbol{def:runtime}},
 	description={\linkonly{def:runtime}},
	sort={T}
}
\newglossaryentry{sym:threshold-residual}{
	type=symbolslist,
 	name={\glsentrysymbol{def:threshold-residual}},
 	description={\linkonly{def:threshold-residual}},
	sort={delta}
}

\newglossaryentry{def:descentdirection}{
	type=definitionlist,
	symbol={\ensuremath{\dual{\Vect{d}}}},
	name={Landweber descent direction},
	description={Landweber method's search direction},
	sort=d
}

\newglossaryentry{def:descentdirection-precursor}{
	type=definitionlist,
	symbol={\ensuremath{\dual{\Vect{R}}}},
	name={Landweber descent direction in $\DualSpace{Y}$},
	description={Landweber method's search direction before applying the dual of the operator $\Mat{A}$},
	sort=Rs
}

\newglossaryentry{def:iterate}{
	type=definitionlist,
	symbol={\ensuremath{\Vect{x}}},
	name={iterate},
	description={iterate generated by the update formula of the iterative method, converge against the true solution},
	sort=x
}
 
\newglossaryentry{def:iterate-dual}{
	type=definitionlist,
	symbol={\ensuremath{\dual{\Vect{x}}}},
	name={dual of the iterate},
	description={iterate projected via the duality mapping associated with the space of the iterate},
	sort=xs
}

\newglossaryentry{def:modulus-smoothness}{
	type=definitionlist,
	symbol={\ensuremath{\rho}},
	name={modulus of smoothness},
	description={it measures the distance of the sum two vectors on the unit sphere from the upper bound $2$},
	sort=rho
}

\newglossaryentry{def:orthogonalization-coefficient}{
	type=definitionlist,
	symbol={\ensuremath{s}},
	name={orthogonalization coefficient},
	description={scalar for the older search directions to subtract from the current search direction},
	sort=s
}

\newglossaryentry{def:orthogonalization-coefficient-bregman}{
	type=definitionlist,
	symbol={\ensuremath{s^{(\text{b})}}},
	name={Bregman orthogonalization coefficient},
	description={scalar for the older search directions to subtract from the current search direction derived via Bregman projections},
	sort=s
}

\newglossaryentry{def:orthogonalization-coefficient-metric}{
	type=definitionlist,
	symbol={\ensuremath{s^{(\text{m})}}},
	name={metric orthogonalization coefficient},
	description={scalar for the older search directions to subtract from the current search direction derived via metric projections},
	sort=s
}

\newglossaryentry{def:residual}{
	type=definitionlist,
	symbol={\ensuremath{R}},
	name={residual},
	description={residual of the inverse problem},
	sort=R
}

\newglossaryentry{def:residualdirection}{
	type=definitionlist,
	symbol={\ensuremath{\Vect{R}}},
	name={residual direction in $\Space{Y}$},
	description={residual direction of the inverse problem},
	sort=R
}

\newglossaryentry{def:searchdirection}{
	type=definitionlist,
	symbol={\ensuremath{\dual{\Vect{u}}}},
	name={search direction in $\DualSpace{X}$},
	description={search direction by which the dual of the current iterate is modified to produce the next iterate},
	sort=u
}

\newglossaryentry{def:searchdirection-orthogonalized}{
	type=definitionlist,
	symbol={\ensuremath{\dual{\Vect{v}}}},
	name={orthogonalized search direction in $\DualSpace{X}$},
	description={search direction made orthogonal in a dual pairing-sense against older search directions via metric projection},
	sort=v
}

\newglossaryentry{def:searchdirection-orthogonalized-precursor}{
	type=definitionlist,
	symbol={\ensuremath{\dual{\Vect{o}}}},
	name={orthogonalized search direction in $\DualSpace{Y}$},
	description={search direction made orthogonal in a dual pairing-sense against older search directions via metric projection before applying the dual of the operator $\Mat{A}$},
	sort=o
}

\newglossaryentry{def:searchdirection-precursor}{
	type=definitionlist,
	symbol={\ensuremath{\dual{\Vect{w}}}},
	name={search direction in $\DualSpace{Y}$},
	description={search direction before applying the dual of the operator $\Mat{A}$},
	sort=w
}

\newglossaryentry{def:stepwidth}{
	type=definitionlist,
	symbol={\ensuremath{\mu}},
	name={stepwidth},
	description={scalar for the search direction to modify the iterate with},
	sort=mu
}

\newglossaryentry{def:truesolution}{
	type=definitionlist,
	symbol={\ensuremath{\widehat{\Vect{x}}}},
	name={true solution},
	description={true solution of the inverse problem $\Mat{A}\Vect{x} = \Vect{y}$},
	sort=xhat
}

\newglossaryentry{sym:descentdirection}{
	type=symbolslist,
	name={\glsentrysymbol{def:descentdirection}},
	description={\linkonly{def:descentdirection}},
	sort={d}
}

\newglossaryentry{sym:descentdirection-precursor}{
	type=symbolslist,
	name={\glsentrysymbol{def:descentdirection-precursor}},
	description={\linkonly{def:descentdirection-precursor}},
	sort={e}
}

\newglossaryentry{sym:iterate}{
	type=symbolslist,
	name={\glsentrysymbol{def:iterate}},
	description={\linkonly{def:iterate}},
	sort={x}
}

\newglossaryentry{sym:iterate-dual}{
	type=symbolslist,
	name={\glsentrysymbol{def:iterate-dual}},
	description={\linkonly{def:iterate-dual}},
	sort={xs}
}

\newglossaryentry{sym:modulus-smoothness}{
	type=symbolslist,
	name={\glsentrysymbol{def:modulus-smoothness}},
	description={\linkonly{def:modulus-smoothness}},
	sort={rho}
}

\newglossaryentry{sym:orthogonalization-coefficient}{
	type=symbolslist,
	name={\glsentrysymbol{def:orthogonalization-coefficient}},
	description={\linkonly{def:orthogonalization-coefficient}},
	sort={s}
}

\newglossaryentry{sym:orthogonalization-coefficient-bregman}{
	type=symbolslist,
	name={\glsentrysymbol{def:orthogonalization-coefficient-bregman}},
	description={\linkonly{def:orthogonalization-coefficient-bregman}},
	sort={sb}
}

\newglossaryentry{sym:orthogonalization-coefficient-metric}{
	type=symbolslist,
	name={\glsentrysymbol{def:orthogonalization-coefficient-metric}},
	description={\linkonly{def:orthogonalization-coefficient-metric}},
	sort={sm}
}

\newglossaryentry{sym:residual}{
	type=symbolslist,
	name={\glsentrysymbol{def:residual}},
	description={\linkonly{def:residual}},
	sort={R}
}

\newglossaryentry{sym:searchdirection}{
	type=symbolslist,
	name={\glsentrysymbol{def:searchdirection}},
	description={\linkonly{def:searchdirection}},
	sort={u}
}

\newglossaryentry{sym:searchdirection-orthogonalized}{
	type=symbolslist,
	name={\glsentrysymbol{def:searchdirection-orthogonalized}},
	description={\linkonly{def:searchdirection-orthogonalized}},
	sort={o}
}

\newglossaryentry{sym:searchdirection-orthogonalized-precursor}{
	type=symbolslist,
	name={\glsentrysymbol{def:searchdirection-orthogonalized-precursor}},
	description={\linkonly{def:searchdirection-orthogonalized-precursor}},
	sort={v}
}

\newglossaryentry{sym:searchdirection-precursor}{
	type=symbolslist,
	name={\glsentrysymbol{def:searchdirection-precursor}},
	description={\linkonly{def:searchdirection-precursor}},
	sort={w}
}

\newglossaryentry{sym:stepwidth}{
	type=symbolslist,
	name={\glsentrysymbol{def:stepwidth}},
	description={\linkonly{def:stepwidth}},
	sort={mu}
}

\newglossaryentry{sym:truesolution}{
	type=symbolslist,
	name={\glsentrysymbol{def:truesolution}},
	description={\linkonly{def:truesolution}},
	sort={xhat}
}

\newcommand{\mypath}{./experiments}%

\title{A \glsentrytext{glos:CG}-type method in Banach spaces with an application to computerized tomography}

\author{Frederik Heber\footnote{Department of Mathematics, Saarland University, 66123 Saarbr\"ucken, Germany ({\tt Email: heber@math.uni-sb.de})}, Frank Sch\"opfer\footnote{Department of Mathematics, University of Oldenburg, 26129 Oldenburg, Germany ({\tt Email: frank.schoepfer@uni-oldenburg.de})}, and Thomas Schuster\footnote{Department of Mathematics, Saarland University, 66123 Saarbr\"ucken, Germany ({\tt Email: thomas.schuster@num.uni-sb.de}), correspondent author}}

\begin{document}

\maketitle

\begin{abstract}
\glostext{CG} methods are one of the most effective iterative methods to solve linear equations in Hilbert spaces. So far, they have been inherently bound to these spaces since they make use of the inner product structure. In more general Banach spaces one of the most prominent iterative solvers are Landweber-type methods that essentially resemble the \glosname{SD} method applied to the normal equation. More advanced are subspace methods that take up the idea of a Krylov-type search space, wherein an optimal solution is sought. However, they do not share the conjugacy property with \glostext{CG} methods. 
In this article we propose that the \glostext{SESOP} method can be considered as an extension of \glostext{CG} methods to Banach spaces. We employ metric projections to orthogonalize the current search direction with respect to the search space from the last iteration. For the $\ell_2$-space our method then exactly coincides with the Polak-Ribi\`ere type of the \glostext{CG} method when applied to the normal equation. We show that such an orthogonalized search space still leads to weak convergence of the subspace method. Moreover, numerical experiments on a random matrix toy problem and 2D computerized tomography on $\ell_{\sym{norm-X}}$-spaces show superior convergence properties over all $p$ compared to non-orthogonalized search spaces. This especially holds for $\ell_{\sym{norm-X}}$-spaces with small $p$. We see that the closer we are to an $\ell_2$-space, the more we recover of the conjugacy property that holds in these spaces, i.\,e., as expected, the more the convergence behaves independently of the size of the truncated search space.
\end{abstract}

\section{Introduction}\label{sec:introduction}

We consider two Banach spaces $\Space{X}$ and $\Space{Y}$ with a continuous linear operator
\begin{equation}
	\label{math:operator-linear}
	\Mat{A} : \Space{X} \rightarrow \Space{Y},
\end{equation}
where the goal is to (iteratively) solve the inverse problem
\begin{equation}
	\label{math:inverse-problem}
	\Mat{A} \Vect{x} = \Vect{y}.
\end{equation}
$\Space{X}$ is assumed to be smooth and uniformly convex and hence $\Space{X}$ is reflexive and has a strictly convex and uniformly smooth dual $\DualSpace{X}$. $\Space{Y}$ can be arbitrary.
The problem~\eqref{math:inverse-problem} may be ill-posed and thus not suitable for direct inversion of the operator. Hence, a regularization scheme is required to obtain a stable solution.

In \cite{Schopfer2006} the Landweber method, well-known and thoroughly investigated in Hilbert spaces, see references therein, has been extended to this setting. It is essentially a steepest descent method on the normal equation
\begin{equation}
	\label{math:normal-equation}
	\dual{\Mat{A}} \Mat{A} \Vect{x} = \dual{\Mat{A}}\Vect{y}
\end{equation}
that results in the minimum-norm solution. However, the method usually suffers from tremendously slow convergence.

Therefore, in \cite{Schopfer2008} the idea of projection onto subspaces has been taken from the family of \glostext{CG} methods. There, multiple simultaneous search directions are admitted to speed up the convergence. These are the so-called subspace methods that differ by the specific choice of search directions. A convex, differentiable line search functional is derived by considering Bregman projections onto hyperplanes that can be efficiently minimized by standard optimization methods, see \cite{Nocedal1999}. This is the \glosfirst{SESOP} method, proposed by \cite{Narkiss2005}, where also weak convergence has been shown. Later, by taking specific search directions into account also strong convergence has been proven, see \cite[Prop.~1]{Schopfer2009}.

However, the set of search directions is still not optimal. In the thesis of \cite{Schopfer2007} the search directions are further modified in the notion of maximizing their pair-wise orthogonality. 
In this article we give its detailed derivation, which also connects the \glostext{SESOP} method with the \glosfirst{CG} family of methods in Hilbert spaces. Moreover, we propose that the subspace methods also generalize the \glostext{CG} methods from Hilbert spaces to Banach spaces. 

In the following we first review some general properties of duality mappings, uniformly smooth Banach spaces, and Bregman distances. 
Next, we repeat the \glostext{SESOP} method and introduce semi-orthogonalized subspaces. 
Then, we prove weak convergence of the method with these semi-orthogonalized search space consisting of the last $\sym{number-search-directions}$ Landweber descent directions. 
We conclude with numerical experiments on a random matrix toy problem and 2D computerized tomography to exemplify the significantly improved convergence compared to \glostext{SESOP}~\cite{Schopfer2008} with no orthogonalized search spaces.

\section{Preliminaries}\label{sec:preliminaries}

Throughout the paper let $\Space{X}$ and $\Space{Y}$ be real Banach spaces with duals $\DualSpace{X}$ and $\DualSpace{Y}$. The space $\Space{X}$ shall be smooth and uniformly convex with a sequentially weak-to-weak continuous duality mapping. The space $\Space{Y}$ is on the other hand arbitrary. Their norms will be denoted by $\norm{.}_{\Space{X}}$ and $\norm{.}_{\DualSpace{X}}$, respectively. 
For $\Vect{x} \in \Space{X}$ and $\dual{\Vect{x}} \in \DualSpace{X}$, we write $\skp{ \Vect{x} }{ \dual{\Vect{x}} } = \skp{ \dual{\Vect{x}} }{ \Vect{x} } = \dual{\Vect{x}}(\Vect{x})$. By $\Space{L}(\Space{X},\Space{Y})$ we denote the space of all continuous linear operators $\Mat{A} : \Space{X} \rightarrow \Space{Y}$ and write $\dual{\Mat{A}}$ for its dual operator $\dual{\Mat{A}} \in \Space{L}(\DualSpace{Y}, \DualSpace{X})$ and $\norm{\Mat{A}} = \norm{\dual{\Mat{A}}}$ for the operator norm of $\Mat{A}$.\footnote{In a finite-dimensional setting, $\Mat{A}$ can be represented as a matrix with its adjoint given by the transposed matrix.}

For real numbers $a$, $b$, we write
\begin{equation*}
	a \vee b = \max \{a,b\}, \quad a \wedge b = \min \{a,b\}.
\end{equation*}
Also, let $p, \dual{p}, r, {\dual{r}} \in (1,\infty)$ be conjugate exponents so that
\begin{equation*}
	\frac 1 p + \frac 1 {\sym{norm-dualX}} = 1 \quad \text{and} \quad \frac 1 r + \frac 1 {\dual{r}} = 1.
\end{equation*}
Let us further note by $p$ and $\sym{norm-dualX}$ the conjugate exponents satisfying
\begin{equation}
	\label{math:conjugated-values-equivalences}
	\tfrac 1 p + \tfrac 1 {\sym{norm-dualX}} = 1 \Leftrightarrow 1 - \tfrac 1 p = \tfrac {p-1} p = \tfrac 1 {\sym{norm-dualX}} \Leftrightarrow \tfrac p {p-1} = {\sym{norm-dualX}} \Leftrightarrow \tfrac p {\sym{norm-dualX}}= p-1 \Leftrightarrow p + {\sym{norm-dualX}} = {\sym{norm-dualX}}p,
\end{equation}
where all of these hold also with $p$ and ${\sym{norm-dualX}}$ exchanged due to symmetry.

Let us briefly motivate the need for the definitions to follow and also explain the general iterative scheme: 
We want to find an element $\Vect{x}$ in a space $\Space{X}$ fulfilling \eqref{math:inverse-problem}. In order to assess distances between elements we require a metric or simpler a norm. We will look at a sequence $\{\sym{iterate}_n\}_n$ of iterates $\sym{iterate}_n$ in the space that will converge against the true solution $\sym{truesolution}$. We desire this fixed point $\sym{truesolution}$ of the iteration to be contained in that space, i.\,e.~it must be closed with respect to the norm. Hence, the space $\Space{X}$ we are looking at is a \emph{Banach space}. 

Convergence will be proven with respect to the \emph{Bregman distance} and not with respect to the norm, where proof attempts have failed so far. The Bregman distance measures the non-linearity gap between two elements.
A further essential ingredient for a proof of convergence is a geometrical inequality that relates the distance between two vectors $\norm{\Vect{x} - \Vect{y}}$ and the norm of each vector $\norm{\Vect{x}}$ and $\norm{\Vect{y}}$, respectively. These are given by the \emph{Xu-Roach inequalities} in uniformly smooth Banach spaces.

We will provide an update formula that yields the next iterate $\sym{iterate}_{n+1}$ in the sequence $\{\sym{iterate}_n\}_n$ given the current one $\sym{iterate}_n$. This update formula requires an additional search direction used to modify the current iterate.

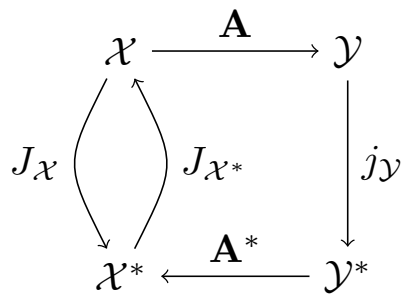
\begin{figure}[htbp]
	\begin{center}
	\resizebox{0.5\textwidth}{!}{
		\ifdefined\DualSpace{\empty}\else\def\DualSpace{\empty}\fi%
\ifdefined\Mat{\empty}\else\def\Mat{\empty}\fi%
\ifdefined\Space{\empty}\else\def\Space{\empty}\fi%
\ifdefined\dual{\empty}\else\def\dual{\empty}\fi%

\begin{tikzpicture}[node distance=2cm]
\node (X) {$\Space{X}$};	
\node[right of=X] (Y) {$\Space{Y}$};	
\node[below of=X] (Xstar) {$\DualSpace{X}$};
\node[below of=Y] (Ystar) {$\DualSpace{Y}$};
\draw[->] (X) -- (Y) node[midway,anchor=south] (A) {$\Mat{A}$};
\draw[->] (Ystar) -- (Xstar) node[midway,anchor=south] (Astar) {$\dual{\Mat{A}}$};
\draw[->] (X) .. controls (-.5,-1) .. (Xstar) node[midway,anchor=east] (Jx) {$J_\Space{X}$};
\draw[->] (Xstar) .. controls (.5,-1) .. (X) node[midway,anchor=west] (Jx) {$J_{\DualSpace{X}}$};
\draw[->] (Y) -- (Ystar) node[midway,anchor=west] (j) {${j_\Space{Y}}$};
\end{tikzpicture}
	}
	\end{center}
	\caption{Relations between the Banach spaces and dual spaces involved in the inverse problem $\Mat{A}\Vect{x}=\Vect{y}$.}
	\label{fig:space-relations}
\end{figure}

The residual $(\Mat{A} \sym{iterate}_n - \Vect{y})$ is a natural choice for such a search direction as it states how much we are off the true solution $\sym{truesolution}$ in the operator's range: $\Mat{A} \sym{iterate}_n - \Vect{y} = \Mat{A} \left ( \sym{iterate}_n - \sym{truesolution} \right )$. As we are only given $\Vect{y}$ in the problem setting, this is all we can do. The residual, however, cannot be used right away as a gradient direction as it resides in $\Space{Y}$ and not in $\Space{X}$.
Now, if the Banach space is not a Hilbert space, then the space and its dual are not isometrically isomorph to one another. In such a case, if the operator $\Mat{A}$ maps an element $\Vect{x} \in \Space{X}$ to an element $\Vect{y}$ in a different space $\Space{Y}$, then its dual operator $\dual{\Mat{A}}$ will not map back in the same spaces but between their respective dual spaces, see Figure~\ref{fig:space-relations}.

To this end, we need operators mapping from a space to its dual and back. Furthermore, these operators should be connected to the norms of the spaces, as the minimum-norm requirement ensures us uniqueness of the solution, see \eqref{math:normal-equation}. This role is fulfilled by so-called \emph{duality mappings}, see Figure~\ref{fig:space-relations} that can be considered a generalization of a norm's gradient, namely the \emph{subgradient}, to the case when it is not differentiable.

Finally, for a gradient method a step width must be known that scales the current search direction. The Landweber method uses a fixed stepwidth but several dynamic, mostly residual-dependent widths are known as well. Naturally, these are either very expensive to calculate or not optimal. Alternatively, for the subspace methods a line search functional is minimized, where even multiple search directions are admitted. This functional consists of computing a \emph{\mbox{Bregman} projection} onto the intersection of planes defined by the search directions. While multiple search directions improve convergence, they are not related in any special way to one another. 

This is then where this article continues: we propose to further use metric projections to semi-\emph{orthogonalize} search directions to enhance the convergence's effectiveness and to keep a memory of outdated search directions. This essentially is the conjugacy property of the \glostext{CG} methods.

	\subsection{Duality Mappings}\label{sec:preliminaries-duality-mappings}

	We recall the definition of duality mappings and some of their properties, all of which can be found in the comprehensive book \cite{Cioranescu1990}.
	
	\begin{definition}[Duality Mapping]\label{def:duality-mapping}
		The mapping $\sym{mapping-duality}_{\sym{norm-X}} : \Space{X} \rightarrow 2^{\DualSpace{X}}$ defined by
		\begin{equation}
			\sym{mapping-duality}_{\sym{norm-X}}(\Vect{x}) = \{ \dual{\Vect{x}} \in \DualSpace{X} \mid \skp{ \dual{\Vect{x}} }{ \Vect{x} }= \norm{ \Vect{x} }^{\sym{norm-X}}, \norm{ \dual{\Vect{x}} } = \norm{ \Vect{x} }^{p-1} \}
		\end{equation}
		is the \emph{duality mapping} of $\Space{X}$ with gauge function $t \mapsto t^{p-1}$.
	\end{definition}

	We remark that the sets $\sym{mapping-duality}_{\sym{norm-X}}(\Vect{x})$ are non-empty by the Hahn-Banach theorem and that they are in general set-valued. We refer to \cite{Schopfer2008} for examples for Hilbert, $L_{\sym{norm-X}}$-spaces, and $\ell_{\sym{norm-X}}$-spaces. By the theorem of Asplund, see \cite[Thm.~4.4]{Cioranescu1990}, duality mappings are the subdifferential of the norm in the form
	\begin{equation}\label{math:mapping-duality-definition}
		\sym{mapping-duality}_{\sym{norm-X}}(\Vect{x}) = \partial \left ( \tfrac 1 p \norm{ \Vect{x} }^{\sym{norm-X}} \right ),
	\end{equation}
	where the subdifferential is the set of subgradients $\dual{\Vect{x}}$ of a function $f$ that each fulfil
	\begin{equation}
		f(\Vect{y}) - f(\Vect{x}) \geq \skp{ \dual{\Vect{x}} }{ \Vect{y} - \Vect{x} } \quad \text{for all } \Vect{y} \in \Space{X}.
	\end{equation}
	That is the theorem of Asplund associates the norm's (sub)gradient in the form \eqref{math:mapping-duality-definition} with the duality mapping of its space.

\pagebreak[0]	\begin{proposition}[Single-valued duality mapping, {\cite[Thm.~3.5, 4.5]{Cioranescu1990}}]\label{prop:smoothness-duality-mapping}
		$\Space{X}$ is smooth if and only if its norm is G{\^a}teaux differentiable on $\Space{X}\setminus 0$ iff $\sym{mapping-duality}_{\sym{norm-X}}$ is single-valued.
	\end{proposition}
	
	\subsection{Uniform Smoothness}\label{sec:preliminaries-uniform-smoothness}
	
	 We will show the convergence of \glostext{SESOP} with orthogonalized search directions in uniformly smooth spaces. 
	The definition of uniform smoothness hinges on the the modulus of smoothness.
	\pagebreak[0]\begin{definition}[{{Uniformly smooth, \cite[Def.~II.1.e.1]{Lindenstrauss1963}}}]\label{def:uniformly-smooth}
		\leavevmode
		Let $\Space{X}$ be a Banach space with $\dim{\Space{X}}\geq 2$.
		\begin{enumerate}[label*=(\alph*)]
			\item\label{itm:modulus-smoothness} Its \emph{modulus of smoothness} is defined by
			\begin{align}
				\label{math:modulus-of-smoothness}
				\sym{modulus-smoothness}_{\Space{X}} (\tau) = \tfrac 1 2 \sup_{\norm{\Vect{x}}=1, \norm{\Vect{y}}=1} \left ( \norm{\Vect{x} + \tau \Vect{y}} + \norm{\Vect{x} - \tau \Vect{y}} - 2 \right ) \quad \tau > 0.
			\end{align}
			\item\label{itm:uniformly-smooth} $\Space{X}$ is said to be uniformly smooth if
			\begin{equation}
				\lim_{\tau \rightarrow \infty} \frac{ \sym{modulus-smoothness}_{\Space{X}} (\tau) }{\tau} = 0.
			\end{equation}
		\end{enumerate}
	\end{definition}
	
	Note that $L_{\sym{norm-X}}$-spaces and $\ell_{\sym{norm-X}}$-spaces with $\sym{norm-X} \in (1,\infty)$ are uniformly convex, see \cite{Clarkson1936}, and thereby also uniformly smooth, see \cite[Lemma~6.7]{Day1944}, see also \cite[\citepage~63]{Lindenstrauss1979}. On the other hand, $\ell_1$ and $l_{\infty}$ are not reflexive\footnote{Note that finite-dimensional spaces are always reflexive.}, see \cite[\citesection~VIII.\,5]{Schechter1971}, and hence can be neither uniformly smooth nor uniformly convex.


	
	In order to prove convergence we will heavily rely on the geometrical characteristics of Banach spaces. Essentially, we need to estimate the distance in the norm between two vectors $\norm{ \Vect{x} - \Vect{y} }$ from both norms $\norm{\Vect{x}}, \norm{\Vect{y}}$ alone. These characteristics are captured in the Xu-Roach inequalities, one of which we repeat here for convenience from \cite{Schopfer2007}.
	\pagebreak[0]\begin{theorem}[{{\cite[Theorem~2, Remark~4]{Xu1991}}}]\label{thm:xu-roach-inequality-smoothness}
		If $\Space{X}$ is uniformly smooth, then for all $\Vect{x}, \Vect{y} \in \Space{X}$, we have
		\begin{equation}\label{math:xu-roach-inequality-smoothness}
			\norm{ \Vect{x} - \Vect{y} }^{\sym{norm-X}} \leq \norm{\Vect{x}}^{\sym{norm-X}} - p \skp{ \sym{mapping-duality}_{\sym{norm-X}} (\Vect{x}) }{ \Vect{y} } + \widetilde{\sigma}_{\sym{norm-X}}( \Vect{x}, \Vect{y} )
		\end{equation}
		with 
		\begin{equation}\label{math:sigma-smoothness}
			\widetilde{\sigma}_{\sym{norm-X}}( \Vect{x}, \Vect{y} ) = \sym{norm-X} G_{\sym{norm-dualX}} \int^1_0 \frac{ (\norm{ \Vect{x} - t \Vect{y} } \vee \norm{ \Vect{x} } )^{\sym{norm-X}} }{ t } \sym{modulus-smoothness}_{\Space{X}} \left ( \frac{ t \norm{\Vect{y}} }{ \norm{ \Vect{x} - t \Vect{y} } \vee \norm{ \Vect{x} } } \right ) dt
		\end{equation}
		where 
		\begin{align*}
			G_{\sym{norm-dualX}} &= 8 \vee 64cK_{\sym{norm-dualX}}^{-1}, \\
			K_{\sym{norm-dualX}} &= 4(2+\sqrt{3})\min\{ \tfrac 1 2 {\sym{norm-dualX}}({\sym{norm-dualX}}-1) \wedge 1, ( \tfrac 1 2 {\sym{norm-dualX}} \wedge 1 )({\sym{norm-dualX}}-1), \\
			& ({\sym{norm-dualX}}-1)\left [ 1-(\sqrt{3} - 1)^{\sym{norm-X}} \right ], 1 - \left [ 1 + (2-\sqrt{3})\sym{norm-X} \right ]^{1-{\sym{norm-dualX}}} \},
		\end{align*}
		and
		\begin{equation}
			c = \frac{ 4\tau_0 }{ \sqrt{ 1 + \tau_0^2 } -1} \prod_{j=1}^{\infty} \left ( 1+ \frac{15\tau_0}{2^{j+2} } \right )\quad \text{with } \tau_0 = \frac{ \sqrt{339} - 18 }{ 30 }.
		\end{equation}
		\begin{proof}
			The proof is given in \cite{Xu1991}. Note that we have additionally used $\sym{modulus-smoothness}_{\Space{X}} (\tau) \geq \sqrt{1 + \tau^2} -1$ for every Banach space $\Space{X}$, see \cite[\citepage~243]{Lindenstrauss1963}.
		\end{proof}
	\end{theorem}
	
	 We also give here a technical lemma on an upper bound on the function $\widetilde{\sigma}_{\sym{norm-X}}$ here, which we require later in the convergence proof.
	\pagebreak[0]\begin{lemma}[Upper bound on $\widetilde{\sigma}$,{\cite[Proof of Prop.~2.39]{Schopfer2007}}]\label{lem:sigma-upperbound}
		Let $\DualSpace{X}$ be a uniformly smooth Banach space with duality mapping $\sym{mapping-duality}_{\sym{norm-dualX}}$. 
		If $0 \neq \Vect{x} \in \Space{X}$, $0 \neq \Mat{A} \in {\cal L}(\Space{X},\Space{Y})$ and $0 \neq \dual{\Vect{y}} \in \DualSpace{Y}$ with an arbitrary Banach space $\Space{Y}$ are given and $\mu > 0$ is defined by
		\begin{equation}\label{math:stepwidth-choice-landweber}
			\mu := \frac{ \tau }{  \norm{\Mat{A}} } \frac{ \norm{ \Vect{x} }^{\sym{norm-X}-1} }{ \norm{ \dual{\Vect{y}} } } \quad \text{for some }\tau \in (0, 1],
		\end{equation}
		then the following estimate is valid:
		\begin{equation}
			\tfrac 1 {\sym{norm-dualX}} \widetilde{\sigma}_{\sym{norm-dualX}} ( \sym{mapping-duality}_{\sym{norm-dualX}} (\Vect{x}), \mu \dual{\Mat{A}} \dual{\Vect{y}}) \leq 2^{\sym{norm-dualX}} G_{\sym{norm-X}} \norm{ \Vect{x} }^{\sym{norm-X}} \sym{modulus-smoothness}_{\DualSpace{X}} (\tau),
		\end{equation}
		whereby $G_{\sym{norm-dualX}}$ is the constant appearing in \eqref{math:sigma-smoothness} and $\sym{modulus-smoothness}_{\DualSpace{X}}$ is the modulus of smoothness of $\DualSpace{X}$, see \cite{Cioranescu1990} for definitions and details.
		\begin{proof}
			We look at the definition of \eqref{math:sigma-smoothness} for $\widetilde{\sigma}_{\sym{norm-dualX}}( \sym{mapping-duality}_{\sym{norm-dualX}} (\Vect{x}), \dual{\Mat{A}} \dual{\Vect{y}} )$ and estimate term by term.
			By \eqref{math:stepwidth-choice-landweber} we get
			\begin{equation*}
				\norm{ \sym{mapping-duality}_{\sym{norm-dualX}} (\Vect{x}) - t \mu \dual{\Mat{A}} \dual{\Vect{y}} } \leq \norm{ \Vect{x} }^{\sym{norm-X}-1} + \mu \norm{ \Mat{A} } \norm{ \dual{\Vect{y}} } \leq 2 \norm{ \Vect{x} }^{\sym{norm-X}-1}
			\end{equation*}
			and
			\begin{equation*}
				\norm{ \sym{mapping-duality}_{\sym{norm-dualX}} (\Vect{x}) - t \mu \dual{\Mat{A}} \dual{\Vect{y}} } \vee \norm{ \sym{mapping-duality}_{\sym{norm-dualX}} (\Vect{x}) } \begin{cases} \leq 2 \norm{ \Vect{x} }^{\sym{norm-X}-1} \\ \geq \norm{ \Vect{x} }^{\sym{norm-X}-1}\end{cases}.
			\end{equation*}
			As the modulus of smoothness $\sym{modulus-smoothness}_{\DualSpace{X}}$ is non-decreasing, see \cite[Prop.~1.e.5]{Lindenstrauss1979}, and with \eqref{math:mapping-duality-definition} and \eqref{math:stepwidth-choice-landweber}, we see that
			\begin{equation*}
				\sym{modulus-smoothness}_{\DualSpace{X}} \left ( \frac{ t \mu \norm{ \dual{\Mat{A}} \dual{\Vect{y}} } }{ \norm{ \sym{mapping-duality}_{\sym{norm-dualX}} (\Vect{x}) - \mu \dual{\Mat{A}} \dual{\Vect{y}} } \vee \norm{ \sym{mapping-duality}_{\sym{norm-dualX}} (\Vect{x}) }  } \right ) \leq \sym{modulus-smoothness}_{\DualSpace{X}} \left ( \frac{ t \mu \norm{ \dual{\Mat{A}} \dual{\Vect{y}} } }{ \norm{ \Vect{x} }^{p-1}  } \right ) \leq \sym{modulus-smoothness}_{\DualSpace{X}} ( t \tau).
			\end{equation*}
			And we finally arrive at the desired estimate,
			\begin{align*}
				\tfrac 1 {\sym{norm-dualX}} \widetilde{\sigma}_{\sym{norm-dualX}} ( \sym{mapping-duality}_{\sym{norm-dualX}} (\Vect{x}), \mu \dual{\Mat{A}} \dual{\Vect{y}}) &\leq 2^{\sym{norm-dualX}} G_{\sym{norm-X}} \norm{ \Vect{x} }^{\sym{norm-X}} \int^1_0 \frac{ \sym{modulus-smoothness}_{\DualSpace{X}} (t\tau) }{ t } dt \\
				&= 2^{\sym{norm-dualX}} G_{\sym{norm-X}} \norm{ \Vect{x} }^{\sym{norm-X}} \int^\tau_0 \frac{ \sym{modulus-smoothness}_{\DualSpace{X}} (t) }{ t } dt \\
				&\leq 2^{\sym{norm-dualX}} G_{\sym{norm-X}} \norm{ \Vect{x} }^{\sym{norm-X}} \sym{modulus-smoothness}_{\DualSpace{X}} (\tau)
			\end{align*}
			as the function $\tau \rightarrow \tfrac{ \sym{modulus-smoothness}_{\DualSpace{X}} (\tau) }{ \tau }$ is non-decreasing, see \cite[Cor.~2.8]{Chidume2009a}
		\end{proof}
	\end{lemma}

	\subsection{Metric and Bregman Projections}\label{sec:preliminaries-projections}

	We summarize essential properties of metric and Bregman projections. For a projection $P$ onto a closed and convex set $\Space{C}$ we have $P^2(\Vect{x}) = P(\Vect{x})$ and $P(\Vect{x}) = \Vect{x} \Leftrightarrow \Vect{x} \in \Space{C}$.

	\pagebreak[0]\begin{definition}[Metric Projection]\label{def:projection-metric}
		The \emph{metric projection} $\sym{projection-metric}$ of $\Vect{x} \in \Space{X}$ onto $\Space{C}$ is the unique element $\sym{projection-metric}_{\Space{C}}(\Vect{x}) \in \Space{C}$ such that
		\begin{equation}
			\norm{ \Vect{x} - \sym{projection-metric}_{\Space{C}}(\Vect{x}) } = \min_{\Vect{y} \in \Space{C}} \norm{ \Vect{x} - \Vect{y} }.	
		\end{equation}
	\end{definition}
	
	Let us then remind of the Bregman distance in the context of generalized distance functions, see also  \cite[\citesection~2.1]{Censor1997}
	\pagebreak[0]\begin{definition}[Bregman Distance]\label{def:distance-bregman}
		For a G{\^a}teaux differentiable convex function $f : \Space{X} \rightarrow \reel$ the function 
		\begin{equation}\label{math:distance-bregman}
			\sym{distance-bregman}_f (\Vect{x}, \Vect{y}) := f(\Vect{y}) - f(\Vect{x}) - \skp{ f'(\Vect{x}) }{ \Vect{y} - \Vect{x} }, \quad \Vect{x}, \Vect{y} \in \Space{X}
		\end{equation}
		is called the \emph{Bregman distance} of $\Vect{x}$ to $\Vect{y}$ with respect to the function $f$.
	\end{definition}

	Here, we consider Bregman distances of functions $f_{\sym{norm-X}}(\Vect{x}) = \tfrac 1 p \norm{ \Vect{x} }^{\sym{norm-X}}$ with $f'_{\sym{norm-X}} = \sym{mapping-duality}_{\sym{norm-X}}$, i.\,e.~the duality mapping of $\Space{X}$, see \eqref{math:mapping-duality-definition}.
	A useful identity~\cite{Schopfer2008} for the Bregman distance is then 
	\begin{equation}
		\label{math:distance-bregman-identity-power}
		\sym{distance-bregman}_{\sym{norm-X}} ( \Vect{x}, \Vect{y} ) = \frac 1 q \norm{ \Vect{x} }^{\sym{norm-X}} - \skp{   \sym{mapping-duality}_{\sym{norm-X}} (\Vect{x}) }{ \Vect{y} } + \frac 1 p \norm{ \Vect{y} }^{\sym{norm-X}}.
	\end{equation}

	\pagebreak[0]\begin{proposition}[Properties of Bregman Distances,{\cite[Theorem~2.12]{Schopfer2006}}]\label{prop:distance-bregman-properties}
		For all \linebreak$\Vect{x}, \Vect{y} \in \Space{X}$ and sequences $\{\Vect{x_n}\}_n$ in $\Space{X}$ the following holds:
		\begin{enumerate}[label*=(\alph*)]
			\item\label{itm:bregman-distance-normlike-properties} $\sym{distance-bregman}_{\sym{norm-X}}(\Vect{x}, \Vect{y}) \geq 0$ and $\sym{distance-bregman}_{\sym{norm-X}}(\Vect{x}, \Vect{y}) = 0 \Leftrightarrow \Vect{x} = \Vect{y}$. 
			\item\label{itm:bregman-distance-boundedness} $\lim_{\norm{ x_n } \rightarrow \infty} \sym{distance-bregman}_{\sym{norm-X}}(\Vect{x_n},\Vect{x}) = \infty$, i.\,e.~the sequence $\{\Vect{x_n}\}_n$ remains bounded if the sequence $\left \{ \sym{distance-bregman}_{\sym{norm-X}}(\Vect{x_n},\Vect{x}) \right \}_n$ is bounded.
			\item\label{itm:bregman-distance-differentiable} $\sym{distance-bregman}_{\sym{norm-X}}$ is continuous in both arguments. It is strictly convex and G{\^a}teaux differentiable with respect to the second variable with $\partial_{\Vect{y}} \sym{distance-bregman}_{\sym{norm-X}}(\Vect{x}, \Vect{y}) = \sym{mapping-duality}_{\sym{norm-X}}(\Vect{y})- \sym{mapping-duality}_{\sym{norm-X}}( \Vect{x})$.
		\end{enumerate}
	\end{proposition}
	
	It is easy to see that in Hilbert spaces metric distance and Bregman distance coincide.

	Finally, equivalent to metric projections in definition~\ref{def:projection-metric}, Bregman projections minimize the Bregman distance with respect to a given convex, non-empty set.
	\pagebreak[0]\begin{definition}[Bregman Projection]\label{def:projection-bregman}
		The \emph{Bregman projection} of $\Vect{x} \in \Space{X}$ onto $\Space{C}$ with respect to the function $f_{\sym{norm-X}}(\Vect{x}) = \frac 1 p \norm{ \Vect{x} }^{\sym{norm-X}}_{\Space{X}}$ is the unique element $\sym{projection-bregman}^{\sym{norm-X}}_{\Space{C}}(\Vect{x}) \in \Space{C}$ such that
		\begin{equation}
		\sym{distance-bregman}_{\sym{norm-X}} \left ( \Vect{x}, \sym{projection-bregman}^{\sym{norm-X}}_{\Space{C}}(\Vect{x}) \right ) = \min_{\Vect{y} \in \Space{X}} \sym{distance-bregman}_{\sym{norm-X}} ( \Vect{x}, \Vect{y}).
		\end{equation}
	\end{definition}

	Finally, let us look at equalities and differences between these two projections.
	\pagebreak[0]\begin{proposition}[{\cite[Prop.~3.6]{Schopfer2008}}]\label{prop:props-bregman-metric-projection}
		\leavevmode
		\begin{enumerate}[label*=(\alph*)]
			\item\label{itm:props-bregman-metric-projection-relation} The Bregman projection and the metric projection are related via
			\begin{equation}
				\sym{projection-metric}_{\Space{C}} (\Vect{x}) - \Vect{x} = \sym{projection-bregman}^{\sym{norm-X}}_{\Space{C}-\Vect{x}}(0) \quad \text{for all }\Vect {x} \in \Space{X}.
			\end{equation}
			Especially we have $\sym{projection-metric}_{\Space{C}} (0) = \sym{projection-bregman}^{\sym{norm-X}}_{\Space{C}} (0)$.
		\item\label{itm:props-bregman-metric-projection-translational} The metric projection has the translational property
		\begin{equation}
			\sym{projection-metric}_{\Vect{y}+\Space{C}}(\Vect{x}) = \Vect{y} + \sym{projection-metric}_{\Space{C}}(\Vect{x}-\Vect{y}) \quad \text{for all }\Vect{x}, \Vect{y} \in \Space{X}.		
		\end{equation}
		This property indeed distinguishes the metric from the Bregman projection since if we had
		$\sym{projection-bregman}^{\sym{norm-X}}_{\Vect{y}+\Space{C}}(x) = \Vect{y} +  \sym{projection-bregman}^{\sym{norm-X}}_{\Space{C}}(\Vect{x}-\Vect{y})$ for all $\Vect{x}, \Vect{y} \in \Space{X}$, then this would already imply their equivalence, $\sym{projection-bregman}^{\sym{norm-X}}_{\Vect{y} + \Space{C}} (\Vect{x}) = \sym{projection-metric}_{\Vect{y}+\Space{C}}(\Vect{x})$ for all $\Vect{x}, \Vect{y} \in \Space{X}$. 
		\end{enumerate}
	\end{proposition}

\section{Methods}\label{sec:methods}

	We now discuss the sequential subspace methods for solving problem~\eqref{math:inverse-problem} given an operator~\eqref{math:operator-linear} and two Banach spaces $\Space{X}$, $\Space{Y}$ and their dual spaces $\DualSpace{X}$, $\DualSpace{Y}$ with properties as stated in the beginning.
	
	
	Let us have first a comment on the notation: In the preliminaries section we stated that space $\Space{X}$ and its dual space $\DualSpace{X}$ are not isometrically isomorph. In consequence, we exploit in the following method the relations between spaces as given in Figure~\ref{fig:space-relations}. To this end, if $\dual{\Vect{v}}$ is an element in the space $\DualSpace{X}$ and can be written as $\dual{\Mat{A}}\dual{\Vect{o}}$, then we call $\dual{\Vect{o}}$ in the following the "precursor" of $\dual{\Vect{v}}$ because of the intimate connection between the spaces $\DualSpace{Y}$ and $\DualSpace{X}$ via the adjoint operator $\dual{\Mat{A}}$. 
	
	Of special importance is the optimality condition, see \cite[Lemma~2.10]{Schopfer2006}, which we repeat here.
	\pagebreak[0]\begin{lemma}[Optimality condition]\label{lem:condition-optimality}
		Let $\Space{X}$ be smooth and uniformly convex and $\Vect{y} \in \Range{\Mat{A}}$.
		\begin{enumerate}[label*=(\alph*)]
			\item\label{itm:minimum-norm-solution-existence} There exists the minimum-norm-solution $\sym{truesolution}$ of \eqref{math:inverse-problem} and $\sym{mapping-duality}_{\sym{norm-X}} (\sym{truesolution}) \in \overline{\Range{\dual{\Mat{A}}}}$.
			\item\label{itm:minimum-norm-solution-optimality} If $\sym{truesolution} \in \Space{X}$ is the minimum-norm-solution of \eqref{math:inverse-problem} and $\widetilde{\Vect{x}} \in \Space{X}$ fulfils $\sym{mapping-duality}_{\sym{norm-X}} (\widetilde{\Vect{x}}) \in \overline{\Range{\dual{\Mat{A}}}}$ and $\sym{truesolution} - \widetilde{\Vect{x}} \in \Null{\Mat{A}}$, then $\widetilde{\Vect{x}} = \sym{truesolution}$.
		\end{enumerate}
	\end{lemma}
		
	\subsection{\glsentryname{glos:SESOP}}\label{sec:methods-SESOP}

	For convenience let us first recall the \glostext{SESOP} method as given in \cite{Schopfer2008} for solving the ill-posed inverse problem $\Mat{A}\Vect{x}=\Vect{y}$ without noise.
	\pagebreak[0]\begin{method}[\glsentrytext{glos:SESOP}]\label{method:sesop}
		\leavevmode
		\begin{enumerate}[label*=(S\arabic*)]
			\item\label{itm:SESOP-start} Take $\sym{iterate}_0$ as initial value with $\sym{mapping-duality}_{\sym{norm-X}} (\sym{iterate}_0) \in \overline{\Range{\dual{\Mat{A}}}}$, set $n := 0$, $\Space{U}_{-1} := \{0\}$ and repeat the following steps:
			\item\label{itm:SESOP-loop-stopcondition} If $\sym{residual}_n := \norm{ \sym{residualdirection}_n } := \norm{ \Mat{A} \sym{iterate}_n - \Vect{y} } = 0$ then STOP else goto \ref{itm:SESOP-loop-update-searchspace}. 
			\item\label{itm:SESOP-loop-update-searchspace} Choose the search space $\Space{U}_n = \Span \{ \sym{searchdirection}_{n,1}, \ldots, \sym{searchdirection}_{n,\sym{number-search-directions}_n} \}  \subset \Range{\dual{\Mat{A}}}$ with $\sym{number-search-directions}_n$ search directions $\sym{searchdirection}_{n,k} \in \Space{U}_n$, $k = 1, \ldots, \sym{number-search-directions}_n$ and with $\sym{number-search-directions}_n$ offsets $\alpha_{n,k} := \skp{ \sym{searchdirection}_{n,k} }{ \Vect{z} }$ for any $\Vect{z} \in \Space{M}_{\Mat{A}\Vect{x}=\Vect{y}} := \{ \Vect{x} \in \Space{X} : \Mat{A} \Vect{x} = \Vect{y} \}$.
			\item\label{itm:SESOP-loop-update-iterate} Compute the new iterate:
			\begin{equation}\label{math:SESOP-update-iterate}
				\sym{iterate}_{n+1} := \sym{mapping-duality}_{\sym{norm-dualX}} \left ( \sym{mapping-duality}_{\sym{norm-X}} (\sym{iterate}_n) - \sum^{\sym{number-search-directions}_n}_{k=1} \sym{stepwidth}_{n,k} \sym{searchdirection}_{n,k} \right )
			\end{equation}
			where $\sym{stepwidth}_n = (\sym{stepwidth}_{n,1}, \ldots, \sym{stepwidth}_{n,\sym{number-search-directions}_n})$ is the solution of the $\sym{number-search-directions}_n$-dimensional optimization problem 
			\begin{equation*}
				\min_{t \in \reel^{\sym{number-search-directions}_n}} h_n(t) 
			\end{equation*}
			with
			\begin{align}\label{math:SESOP-linesearch-functional}
				h_n(t) &:= \frac 1 {\sym{norm-dualX}} \norm{ \sym{mapping-duality}_{\sym{norm-X}} (\sym{iterate}_n) - \sum^{\sym{number-search-directions}_n}_{k=1} t_k \sym{searchdirection}_{n,k} }^{\sym{norm-dualX}} + \sum^{\sym{number-search-directions}_n}_{k=1} t_k \alpha_{n,k} \\
				\label{math:SESOP-linesearch-functional-derivative}
				\partial_j h_n(t) &= - \SKP{ \sym{searchdirection}_{j,k} }{ \sym{mapping-duality}_{\sym{norm-dualX}} \left ( \sym{mapping-duality}_{\sym{norm-X}} (\sym{iterate}_n) - \sum^{\sym{number-search-directions}_n}_{k=1} t_k \sym{searchdirection}_{n,k} \right ) } + \alpha_{j,k} \quad \forall j =1, \ldots, \sym{number-search-directions}_n
			\end{align}
			\item\label{itm:SESOP-loop-continue} Set $n \leftarrow n+1$ and goto \ref{itm:SESOP-loop-stopcondition}. 
		\end{enumerate}
	\end{method}
	Note that \eqref{math:SESOP-linesearch-functional} is strictly convex, hence $\sym{stepwidth}$ is unique.

	Convergence of the method essentially depends on the choice of the search space $\Space{U}_n$ and associated offsets $\alpha_n$ per iteration step $n$, see step~\ref{itm:SESOP-loop-update-searchspace}. We state a few common choices, taken from \cite{Schopfer2008,Schopfer2009}, using the \emph{Landweber descent direction} $\sym{descentdirection}_n := \dual{\Mat{A}} \sym{descentdirection-precursor}_n = \dual{\Mat{A}} \sym{mapping-duality}_{\sym{norm-Y}} \bigl ( \Mat{A} \sym{iterate}_{n} - \Vect{y} \bigr )$ with precursor $\sym{descentdirection-precursor}_n=\sym{mapping-duality}_{\sym{norm-Y}}\sym{residualdirection}_n$. Here $\sym{mapping-duality}_{\sym{norm-Y}}$ denotes a single-valued selection of the set-valued duality mapping of $Y$.
	\begin{enumerate}[label*=(\alph*)]
		\item\label{itm:searchspace-expanding} Expanding: $\Space{U}_n^{\text{exp}} = \Span \{ \sym{descentdirection}_0, \ldots, \sym{descentdirection}_{n} \}$, $\alpha_{n,k}^{\text{exp}} = \skp{ \sym{descentdirection-precursor}_{n,k} }{ \Vect{y} }$ with dimension $|\Space{U}_n| = n+1$
		\item\label{itm:searchspace-truncated} Truncated: $\Space{U}_n^{\text{trunc}} = \Span \{ \sym{descentdirection}_{n-\sym{number-search-directions}_n+1}, \ldots, \sym{descentdirection}_{n} \}$, $\alpha_{n,k}^{\text{trunc}} = \skp{ \sym{descentdirection-precursor}_{n,k} }{ \Vect{y} }$ with dimension $|\Space{U}_n| = \sym{number-search-directions}_n := N \wedge (n+1)$ for some fixed $N \in \nat$
		\item\label{itm:searchspace-nemirovsky1} Nemirovsky I: $\Space{U}_n^{\text{Nem1}} = \Span \{ \sym{descentdirection}_{n}, J_{\Space{X}} (x_n) - J_{\Space{X}} (x_0) \}$, $\alpha_{n,k}^{\text{Nem1}} = \skp{ \dual{\Vect{v}}_{k,l} }{ \Vect{y} }$ with dimension  $|\Space{U}_n| = 2$
		\item\label{itm:searchspace-nemirovsky2} Nemirovsky II: $\Space{U}_n^{\text{Nem2}} = \Span \{ \sym{descentdirection}_{n}, J_{\Space{X}} (x_n) - J_{\Space{X}} (x_{n-1}) \}$, $\alpha_{n,k}^{\text{Nem2}} = \skp{ \dual{\Vect{v}}_{k,l} }{ \Vect{y} }$ with dimension $|\Space{U}_n| = 2$
	\end{enumerate}
	Note that the Nemirovsky directions of cases \ref{itm:searchspace-nemirovsky1} and \ref{itm:searchspace-nemirovsky2} provide strong convergence, \cite[Prop.~1]{Schopfer2009}, are not considered in the scope of this article.
		Also, for finite-dimensional spaces $\Space{X}$ and $\Space{Y}$ weak and strong convergence coincide~\cite[Thm~4.3]{Schechter1971}.

	In the cases \ref{itm:searchspace-expanding} or \ref{itm:searchspace-truncated} the hyperplane offsets can simply be calculated by $\alpha_{n,k} := \skp{ \sym{searchdirection-precursor}_{n,k} }{ \Vect{y} }$.

	\subsection{Generalized \glsentrytext{glos:CG}}\label{sec:methods-generalized-CG}

	Employing multiple search directions greatly improves convergence as indicated by experiments~\cite[\citesection~5]{Schopfer2008}. However, the search directions are still not related in any particular way among one another. In the following we would like to maximize "distinctiveness" of the search directions $\sym{searchdirection}_{n,k}$ in the truncated search space $\Space{U}_n$, i.\,e.~we orthogonalize them with respect to previous ones contained in $\Space{U}_n$ and to some extent to older search directions not present in $\Space{U}_n$. 
	
	Because of the Banach space structure this orthogonality holds only one way: new descent directions are made orthogonal to ones already contained in the search space but not the other way round, i.\,e.~it is not symmetric.
	
	To make such a (semi-)orthogonalized search space distinguishable from $\Space{U}_n^{\text{trunc}}$, we denote it as $\Space{V}_n^{\text{trunc}}$ in the following. Note that in this section we extend the details connecting \glostext{SESOP} and \glostext{CG} as stated in \cite[Sect.~2.6]{Schopfer2007}.

	Based on this semi-orthogonality, we then want to construct a search space, similar to \glostext{SESOP}'s truncated search space $\Space{U}_n^{\text{trunc}}$, using orthogonalized directions. To this end, we examine search directions derived in the following way from the Landweber descent direction $\sym{descentdirection}_n = \dual{\Mat{A}} \sym{mapping-duality}_{\sym{norm-Y}} \bigl ( \Mat{A} \sym{iterate}_{n} - \Vect{y} \bigr )$ to obtain the semi-orthogonalized search space $\Space{V}_n^{\text{trunc}} := \{\sym{searchdirection-orthogonalized}_{n,k}\}^{\sym{number-search-directions}_{n}}_{k=1}$,
	\begin{align}
		\label{math:orthogonalized-search-directions-olddirections}
		\sym{searchdirection-orthogonalized}_{n,k} &:= \sym{searchdirection-orthogonalized}_{n-1,k+1} \quad \text{for } k=1,\ldots,\sym{number-search-directions}_{n}-1
\\
		\label{math:orthogonalized-search-directions}
		\sym{searchdirection-orthogonalized}_{n,\sym{number-search-directions}_n} &:= \sym{descentdirection}_n - \sum^{\sym{number-search-directions}_{n-1}}_{k=1} \sym{orthogonalization-coefficient}_{n,k} \sym{searchdirection-orthogonalized}_{n-1,k} \\
		\label{math:orthogonalized-search-directions-projections}
		&= \sym{descentdirection}_n - \sym{projection-metric}_{\Space{V}_{n-1}^{\text{trunc}}} \left ( \sym{descentdirection}_n \right ).
	\end{align}
	Here, $\sym{orthogonalization-coefficient}_{n,k}$ is an orthogonalization coefficient obtained from metric projection of the Landweber descent direction $\sym{descentdirection}_n$ onto the search space $\Space{V}_{n-1} = \{\sym{searchdirection-orthogonalized}_{n-1,k}\}^{\sym{number-search-directions}_{n-1}}_{k=1}$. For convenience, let us define,
	\begin{equation}
		\label{math:orthogonalization-coefficient-argmin}
		\sym{orthogonalization-coefficient}_{n} = \left (\sym{orthogonalization-coefficient}_{n,1}, \ldots, \sym{orthogonalization-coefficient}_{n,\sym{number-search-directions}_n} \right ) := \argmin_{\sym{orthogonalization-coefficient} \in \reel^{\sym{number-search-directions}_n}} g_n (\sym{orthogonalization-coefficient})
	\end{equation}
	and
	\begin{align}
		\label{math:orthogonalization-coefficient-metric}
		 g_n (\sym{orthogonalization-coefficient}) &:= \norm{ \sym{descentdirection}_n - \sum_{i=1}^{\sym{number-search-directions}_{n-1}} \sym{orthogonalization-coefficient}_i \sym{searchdirection-orthogonalized}_{n-1,i} }^{\sym{norm-dualX}}_{\DualSpace{X}}. 
	\end{align}
	
	In Hilbert spaces this would become the familiar Gram-Schmidt procedure, see also \cite[\citepage~218]{Muscat2014}.

	We show that all search directions $\sym{searchdirection-orthogonalized}_{n,j}$, $j=1,\ldots,\sym{number-search-directions}_n$, are pairwise semi-orthogonal to in a certain sense. This is the Banach space's counterpart of the conjugacy property in a Hilbert space when treating the normal equation, c.\,f.~\cite[p.~102]{Nocedal1999}.
	\pagebreak[0]\begin{corollary}[Semi-orthogonal search directions]\label{cor:orthogonal-searchdirections}
		We have
		\begin{equation}
			\label{math:conjugacy-searchdirections}
			\skp{ \sym{searchdirection-orthogonalized}_{n,j} }{ \sym{mapping-duality}_{\sym{norm-dualX}} (\sym{searchdirection-orthogonalized}_{n,k}) } = 0 \quad \forall 1 \leq j < k \leq \sym{number-search-directions}_n.
		\end{equation}
		\begin{proof}
			This follows directly from the optimality condition on $g_n(\sym{orthogonalization-coefficient})$ for \eqref{math:orthogonalization-coefficient-argmin},
			\begin{align}
				\label{math:orthogonalization-coefficient-metric-derivative}
				0 = \partial_j g_n (\sym{orthogonalization-coefficient}_{\sym{number-search-directions}_n}) &= - \SKP { \sym{searchdirection-orthogonalized}_{n-1,j} }{ \sym{mapping-duality}_{\sym{norm-dualX}} \left ( \sym{descentdirection}_n - \sum_{k=1}^{\sym{number-search-directions}_{n-1}} \sym{orthogonalization-coefficient}_{n,k} \sym{searchdirection-orthogonalized}_{n-1,k} \right ) } \\
				&= - \skp { \sym{searchdirection-orthogonalized}_{n-1,j} }{ \sym{mapping-duality}_{\sym{norm-dualX}} \left ( \sym{searchdirection-orthogonalized}_{n,\sym{number-search-directions}_n} \right ) } \nonumber
			\end{align}
			for all $j=1,\ldots,\sym{number-search-directions}_{n-1}$.
			As this holds for all $n$ and by the successive construction~\eqref{math:orthogonalized-search-directions-olddirections} of the search spaces $\Space{V}_n^{\text{trunc}}$, we have by induction that all search directions in $\Space{V}_n^{\text{trunc}}$ are orthogonal with respect to one another.
		\end{proof}
	\end{corollary}
	Note that this semi-orthogonality becomes a full orthogonality in the case of a Hilbert space $\Space{X}$ with a duality mapping $\sym{mapping-duality}_{\sym{norm-X}}$ of power type $\sym{norm-X}=2$, as $\sym{mapping-duality}_{\sym{norm-X}}$ can then be identified with the identity.
	For simplicity, we speak in the following only of orthogonal search directions. 

	We would like to stress that the summation in \eqref{math:orthogonalized-search-directions} is over all directions in the search space $\Space{V}^{\text{trunc}}_{n-1}$, including $\sym{searchdirection-orthogonalized}_{n-1,1}$, that is not included in $\Space{V}_n^{\text{trunc}}$, see \eqref{math:orthogonalized-search-directions-olddirections}. Otherwise, the orthogonalization would not change the search space but only modify its spanning vectors. Therefore, we clearly have the case of $\Space{U}_n^{\text{trunc}} \not \subset \Space{V}_n^{\text{trunc}}$ for $n\geq 1$ in general.
	Using the expanding search space $\Space{U}_n^{\text{exp}}$ on the other hand, see case~\ref{itm:searchspace-expanding} in section~\ref{sec:methods-SESOP}, the above orthogonalization of the new search direction would \emph{not} change the iteration, i.\,e. for every $n$ we have $\Space{U}_n^{\text{exp}} = \Space{V}_n^{\text{exp}}$. 
		
	\pagebreak[0]\begin{example}
	In order to highlight the notational equivalence with \glostext{CG} methods in Hilbert spaces, let us consider briefly just a single search direction, i.\,e.~$\Space{V}_n^{\text{trunc}} = \Span\{\sym{searchdirection-orthogonalized}_n\}$, per iteration step $n$, we obtain
	\begin{equation}
		\label{math:orthogonalized-search-directions-single-direction}
		\sym{searchdirection-orthogonalized}_n = \sym{descentdirection}_n - \sym{orthogonalization-coefficient}_{n} \sym{searchdirection-orthogonalized}_{n-1}.	
	\end{equation}
	Here, we have the Landweber descent direction $\sym{descentdirection}_n$ that is modified by the last search direction scaled by the orthogonalization coefficient, i.\,e.~$\Vect{g_n} = \Vect{d_n} - s_n \Vect{g_{n-1}}$ where $\Vect{g_n}$ is the current search direction and $\Vect{d_n}$ is the current gradient direction in the usual notation of the \glostext{CG} methods in Hilbert spaces, c.\,f.~\cite[\citechapter~5]{Nocedal1999} and also \cite{Yuan1995}.
	\end{example}

	With each search direction a hyperplane offset is required that relates this hyperplane to the solution manifold $\Space{M}_{\Mat{A}\Vect{x}=\Vect{y}}$, c.\,f. step~\ref{itm:SESOP-loop-update-searchspace} in method~\ref{method:sesop}. Naturally, the last offset changes this step when we orthogonalize the descent direction $\sym{descentdirection}_n$. Hence, instead of $\alpha_{n,\sym{number-search-directions}_n}$ associated to $\sym{descentdirection}_n$, we require a $\beta_{n,\sym{number-search-directions}_n}$ associated to $\sym{searchdirection-orthogonalized}_n$.
	
	Calculating the new hyperplane offsets $\beta_{n,k}$ to each orthogonalized search direction $\sym{searchdirection-orthogonalized}_{n,k}$ in $\Space{V}_n^{\text{trunc}}$ is then,
	\begin{equation}
		\label{math:orthogonalization-offsets}
		\beta_{n,\sym{number-search-directions}_n} := \skp{ \sym{searchdirection-orthogonalized-precursor}_{n,\sym{number-search-directions}_n} }{ \Vect{y} } = \skp{ \sym{descentdirection-precursor}_n - \sum^{\sym{number-search-directions}_{n-1}}_{k=1} \sym{orthogonalization-coefficient}_{n,k} \sym{searchdirection-orthogonalized-precursor}_{n-1,k} }{ \Vect{y} } = \alpha_{n,\sym{number-search-directions}_n} - \sum^{\sym{number-search-directions}_{n-1}}_{k=1} \sym{orthogonalization-coefficient}_{n,k} \beta_{n-1,k},
	\end{equation}
	where we used 
	\begin{align}
		\sym{searchdirection-orthogonalized}_{n,\sym{number-search-directions}_n} &= \sym{descentdirection}_n - \sum^{\sym{number-search-directions}_{n-1}}_{k=1} \sym{orthogonalization-coefficient}_{n,k} \sym{searchdirection-orthogonalized}_{n-1,k} = \sym{descentdirection}_n - \sum^{\sym{number-search-directions}_{n-1}}_{k=1} \sym{orthogonalization-coefficient}_{n,k} \dual{\Mat{A}} \sym{searchdirection-orthogonalized-precursor}_{n-1,k} \nonumber \\
		\label{math:orthogonalized-searchdirections-precursors}
		&= \dual{\Mat{A}} \left ( \sym{descentdirection-precursor}_n - \sum^{\sym{number-search-directions}_{n-1}}_{k=1} \sym{orthogonalization-coefficient}_{n,k} \sym{searchdirection-orthogonalized-precursor}_{n-1,k} \right ) =: \dual{\Mat{A}} \sym{searchdirection-orthogonalized-precursor}_{n,\sym{number-search-directions}_n},
	\end{align}
	with the precursor $\sym{descentdirection-precursor}_n=\sym{mapping-duality}_{\sym{norm-Y}}\sym{residualdirection}_n$ of the Landweber descent direction $\sym{descentdirection-precursor}_n$, $\sym{descentdirection}_n = \dual{\Mat{A}} \sym{descentdirection-precursor}_n$.

	These orthogonalization coefficients $\sym{orthogonalization-coefficient}_{n,k}$ are calculated by minimizing \eqref{math:orthogonalization-coefficient-metric} with the derivative~\eqref{math:orthogonalization-coefficient-metric-derivative} using standard techniques.
	Note that in the proof of Lemma~\ref{lem:descent-direction-property} we give a good starting value for this line search problem.
	Once we know the coefficients $\sym{orthogonalization-coefficient}_{n,k}$ in \eqref{math:orthogonalized-search-directions}, we can easily evaluate \eqref{math:orthogonalization-offsets}, knowing all other offsets $\beta_{n-1,k}$, $k=1,\ldots, \sym{number-search-directions}_{n-1}$, from previous iterations, 

	\subsection{Proof of Convergence}\label{sec:methods-convergence-proof}

		Now we would like to show that, with the semi-orthogonalized truncated search space $\Space{V}^{\text{trunc}}_n$, method~\ref{method:sesop} still converges weakly to a solution of $\Mat{A}\Vect{x}=\Vect{y}$.
		
		We will assume that $\sym{residual}_n \neq 0$ without loss of generality in the following theorem and in some of the corollaries and lemmata in support of this theorem. This is a valid assumption because if we get $\sym{residual}_n = 0$ at some step $n$, then it holds that $\Mat{A}\sym{iterate}_n = \Vect{y}$ and we are done.
		As a direct consequence of this assumption, it follows that the Landweber descent direction is always non-zero.
		\pagebreak[0]\begin{corollary}[Non-zero descent direction]\label{cor:non-zero-descentdirection}
			If $\sym{residual}_n \neq 0$, we have $\sym{descentdirection}_n \neq 0$.
			\begin{proof}
				Assume the contrary, $\sym{descentdirection}_n = 0$, and let be given a $\Vect{z} \in \Space{M}_{\Mat{A}\Vect{x}=\Vect{y}}$, then we get
				\begin{equation*}
					0 = \skp{ \sym{descentdirection}_n }{ \sym{iterate}_n - \Vect{z} } = \skp{ \sym{descentdirection-precursor}_n }{ \Mat{A}\sym{iterate}_n  -\Mat{A}\Vect{z} } = \norm { \sym{residualdirection}_n }^r = \sym{residual}_n^r,
				\end{equation*}
				which is a contradiction.
			
				Note that by the very same argument we also have $\sym{descentdirection-precursor}_n \neq 0$.
			\end{proof}
		\end{corollary}

		Then, we need to prove something similar to the "expanding subspace" property of the \glostext{CG} methods, c.\,f.~\cite[Theorem~5.2]{Nocedal1999}. 
		
		\pagebreak[0]\begin{corollary}[Truncated Subspace Minimization]\label{cor:expanding-subspace}
			At iteration step $n$ let be given an orthogonalized search space $\Space{V}_n^{\text{trunc}} := \Span\{\sym{searchdirection-orthogonalized}_{n,1}, \ldots, \sym{searchdirection-orthogonalized}_{n,\sym{number-search-directions}_n}\}$ with precursors $\{\sym{searchdirection-orthogonalized-precursor}_{n,i}\}_{i=1}^{\sym{number-search-directions}_n}$.
			\begin{enumerate}[label*=(\alph*)]
			\item\label{item:expanding-subspace-orthogonal-residual} We have that all old search directions still contained in $\Space{V}_n^{\text{trunc}}$ are orthogonal with respect to the dual pairing with current and old residuals,
			\begin{align}
				&\skp{ \sym{searchdirection-orthogonalized-precursor}_{n-1,1} }{ \sym{residualdirection}_{n-\sym{number-search-directions}_n} } = 0\nonumber \\
				&\ldots \nonumber \\
				\label{math:residual-orthogonal-oldsearchdirections}
				&\skp{ \sym{searchdirection-orthogonalized-precursor}_{n-1,1} }{ \sym{residualdirection}_n } = 0, \ldots, \skp{ \sym{searchdirection-orthogonalized-precursor}_{n-1,\sym{number-search-directions}_{n-1}} }{ \sym{residualdirection}_n } = 0.
			\end{align}
			Note that this extends to $\{\sym{searchdirection-orthogonalized-precursor}_{n,1}, \ldots, \sym{searchdirection-orthogonalized-precursor}_{n,\sym{number-search-directions}_{n}-1}\}$ by construction of $\Space{V}^{\text{trunc}}_n$.
			\item\label{item:expanding-subspace-nonzero-search-direction} If we have $\sym{residual}_n \neq 0$, then it also holds that $\sym{searchdirection-orthogonalized}_{n,\sym{number-search-directions}_n} \neq 0$.
			\item\label{item:expanding-subspace-linear-independence} The set of vectors $\{\sym{searchdirection-orthogonalized-precursor}_{n,1}, \ldots, \sym{searchdirection-orthogonalized-precursor}_{n,\sym{number-search-directions}_n}\}$ and $\{\sym{searchdirection-orthogonalized}_{n,1}, \ldots, \sym{searchdirection-orthogonalized}_{n,\sym{number-search-directions}_n}\}$ are each linearly independent. 
			\end{enumerate}
			\begin{proof}
				We first prove \ref{item:expanding-subspace-orthogonal-residual}.			
				Assume $\sym{residual}_n \neq 0$, i.\,e.~$\sym{residualdirection}_n \neq 0$. 
				Let us inspect the optimality condition of the step width functional $h_{n-1}(\sym{stepwidth})$, see \eqref{math:SESOP-linesearch-functional}, at $n-1$ for $j=1,\ldots,\sym{number-search-directions}_{n-1}$,
				\begin{align*}
					0 &= \partial_j h_{n-1} (\sym{stepwidth}_{n-1}) \\
					&= - \skp{ \sym{searchdirection-orthogonalized}_{n-1,j} }{ \sym{mapping-duality}_{\sym{norm-dualX}}\left ( \sym{mapping-duality}_{\sym{norm-X}}(\sym{iterate}_{n-1}) - \sum^{\sym{number-search-directions}_{n-1}}_{k=1} \sym{stepwidth}_{n-1,k} \sym{searchdirection-orthogonalized}_{n-1,k} \right ) } + \beta_{n-1,j} \\
					&= - \skp{ \sym{searchdirection-orthogonalized-precursor}_{n-1,j} }{ \Mat{A}\sym{iterate}_n } + \skp{ \sym{searchdirection-orthogonalized-precursor}_{n-1,j} }{ \Vect{y} } \\
					&= - \skp{ \sym{searchdirection-orthogonalized-precursor}_{n-1,j} }{ \sym{residualdirection}_n },
				\end{align*}
				where we needed \eqref{math:SESOP-update-iterate} and \eqref{math:orthogonalization-offsets}. The statement then follows by \eqref{math:orthogonalized-search-directions-olddirections} and stepping back until $n-\sym{number-search-directions}_n$, knowing that
				\begin{equation*}
					\Space{V}^{\text{trunc}}_n = \{ \sym{searchdirection-orthogonalized}_{n-\sym{number-search-directions}_n-1,\sym{number-search-directions}_{n-\sym{number-search-directions}_n}-1}, \ldots, \sym{searchdirection-orthogonalized}_{n, \sym{number-search-directions}_n} \}.
				\end{equation*}

				Then, we can continue with \ref{item:expanding-subspace-nonzero-search-direction}.  By Corollary~\ref{cor:non-zero-descentdirection} we have $\sym{descentdirection}_n \neq 0$ and thereby with \eqref{math:orthogonalized-search-directions} we have to show that $\sym{searchdirection-orthogonalized}_{n,\sym{number-search-directions}_{n}}$ is not contained in $\Space{V}_{n-1}$. To this end, let be $\Vect{z} \in \Space{M}_{\Mat{A}\Vect{x}=\Vect{y}}$ and thus $\sym{iterate}_n - \Vect{z} \neq 0$. Furthermore, let be given $\lambda_1, \ldots, \lambda_{\sym{number-search-directions}_{n-1}}, \sigma \in \reel$ with
				\begin{equation*}
					\sum^{\sym{number-search-directions}_{n-1}}_{k=1} \lambda_k \sym{searchdirection-orthogonalized}_{n-1,k} + \sigma \sym{descentdirection}_{n} = 0.
				\end{equation*}
				Then by using \ref{item:expanding-subspace-orthogonal-residual},
				\begin{align*}
					0 &= \sum^{\sym{number-search-directions}_{n-1}}_{k=1} \lambda_k \skp{ \sym{searchdirection-orthogonalized}_{n-1,k} }{ \sym{iterate}_n - \Vect{z} } + \sigma \skp{ \sym{descentdirection}_{n} }{ \sym{iterate}_n - \Vect{z} } \\
					&= \sum^{\sym{number-search-directions}_{n-1}}_{k=1} \lambda_k \skp{ \sym{searchdirection-orthogonalized-precursor}_{n-1,k} }{ \sym{residualdirection}_n } + \sigma \skp{ \sym{mapping-duality}_{\sym{norm-Y}} (\sym{residualdirection}_n) }{ \sym{residualdirection}_n } \\
					&=\sigma \sym{residual}_n^r,
				\end{align*}
				we get $\sigma = 0$ by contradiction.

				Next, for \ref{item:expanding-subspace-linear-independence} it suffices to show that $\{\sym{searchdirection-orthogonalized}_{n,1}, \ldots, \sym{searchdirection-orthogonalized}_{n,\sym{number-search-directions}_n}\}$ are linearly independent. Assume again $\sym{residual}_n \neq 0$ and let similarly be given $\Vect{z} \in \Space{M}_{\Mat{A}\Vect{x}=\Vect{y}}$ and $\lambda_1, \ldots, \lambda_{\sym{number-search-directions}_n} \in \reel$ with
				\begin{equation*}
					\sum^{\sym{number-search-directions}_n}_{k=1} \lambda_k \sym{searchdirection-orthogonalized}_{n,k} = 0.
				\end{equation*}
				Using \eqref{math:orthogonalized-search-directions} and \ref{item:expanding-subspace-orthogonal-residual} on $\sym{searchdirection-orthogonalized}_{n,\sym{number-search-directions}_n}$ and looking at
				\begin{align*}
					0 &= \sum^{\sym{number-search-directions}_n}_{k=1} \lambda_k \skp{ \sym{searchdirection-orthogonalized}_{n,k} }{ \sym{iterate}_n - \Vect{z} } = \sum^{\sym{number-search-directions}_n-1}_{k=1} \lambda_k \skp{ \sym{searchdirection-orthogonalized-precursor}_{n,k} }{ \sym{residualdirection}_n } + \lambda_{N_n} \skp{ \sym{searchdirection-orthogonalized-precursor}_{n,\sym{number-search-directions}_n} }{ \sym{residualdirection}_n } \\
					&= \lambda_{N_n} \skp{ \sym{searchdirection-orthogonalized-precursor}_{n,\sym{number-search-directions}_n} }{ \sym{residualdirection}_n } = \lambda_{N_n} \skp{ \sym{descentdirection-precursor}_n }{ \sym{residualdirection}_n } - \sum^{\sym{number-search-directions}_n}_{k=1} \sym{orthogonalization-coefficient}_{n,k}\skp{ \sym{searchdirection-orthogonalized-precursor}_{n-1,k} }{ \sym{residualdirection}_n } \\
					&=\lambda_{N_n} \sym{residual}_n^r,
				\end{align*}
				we get $\lambda_{N_n} = 0$. We continue with
				\begin{equation*}
					0 = \sum^{N_n-1}_{k=0} \lambda_k \skp{ \sym{searchdirection-orthogonalized}_{n,k} }{ \sym{iterate}_{n-1} - \Vect{z} } = \ldots = \lambda_{N_n-1} \sym{residual}_{n-1}^r,
				\end{equation*}
				and by induction we get $\lambda_k = 0$ for all $k=1,\ldots,\sym{number-search-directions}_n$. Hence, also the search directions $\sym{searchdirection-orthogonalized}_{n,k}$ are linearly independent.
				
%
			\end{proof}
		\end{corollary}

		Last but not least, we show that our solution manifold is contained in the intersection of hyperplanes and that iterates and search directions obtained via the update formulas~\eqref{math:SESOP-update-iterate} and \eqref{math:orthogonalized-search-directions} derive from Bregman projections onto the intersection of hyperplanes and search space, respectively. While most of this is not needed in the convergence proof, it is very illustrative for the general procedure.
		
		\pagebreak[0]\begin{lemma}[Intersection of hyperplanes]\label{lem:intersection-hyperplanes}
			\leavevmode
			\begin{enumerate}[label*=(\alph*)]
			\item\label{itm:solution-manifold-in-intersection}
			For the solution manifold $\Space{M}_{\Mat{A}\Vect{x} = \Vect{y}}$ it holds
			\begin{equation}\label{math:solution-manifold-in-intersection}
				\Space{M}_{\Mat{A}\Vect{x} = \Vect{y}} \subset \Space{H}_n := \bigcap_{k=1}^{\sym{number-search-directions}_n} \Space{H}(\sym{searchdirection-orthogonalized}_{n, k}, \beta_{n,k})
			\end{equation}
			with the intersection of the hyperplanes $\Space{H}_n$ of all search directions in $\Space{V}_n^{\text{trunc}}$.
			\item\label{itm:dualiterate-in-range-adjoint-operator}
			Also, we have that $ \sym{mapping-duality}_{\sym{norm-X}} (\sym{iterate}_{n}) - \sym{mapping-duality}_{\sym{norm-X}} (\sym{iterate}_{0}) \in \bigcup_{n} \Space{V}_{n} \subset \overline{\Range{\dual{\Mat{A}}}}$ for all $n$.
			\item\label{itm:iterate-projection-intersection}
			Furthermore, the next iterate $\sym{iterate}_{n+1}$ is the Bregman projection of the current iterate $\sym{iterate}_{n}$ onto the intersection,
			\begin{equation}\label{math:iterate-projection-intersection}
				\sym{iterate}_{n+1} = \sym{projection-bregman}^{\sym{norm-X}}_{\Space{H}_n} (\sym{iterate}_n)
			\end{equation}
			and also,
			\begin{equation}\label{math:dualiterate-projection-intersection}
				\sym{mapping-duality}_{\sym{norm-dualX}} ( \sym{iterate}_{n+1} ) = \sym{projection-bregman}^{\sym{norm-X}}_{\sym{mapping-duality}_{\sym{norm-X}} (\sym{iterate}_n ) + \Space{V}_{n}^{\text{trunc}}} \left ( \sym{mapping-duality}_{\sym{norm-X}} ( \Vect{z} ) \right ) \quad \forall \Vect{z} \in \Space{M}_{\Mat{A}\Vect{x}=\Vect{y}}.
			\end{equation}
			\item\label{itm:searchdirection-projection-intersection}
			We have the search direction $\sym{searchdirection-orthogonalized}_{n,\sym{number-search-directions}_n}$ as the Bregman projection of the Landweber descent direction $\sym{descentdirection}_n$,
			\begin{equation}\label{math:searchdirection-projection-intersection}
				\sym{mapping-duality}_{\sym{norm-dualX}} ( \sym{searchdirection-orthogonalized}_{n,\sym{number-search-directions}_n} ) = \sym{projection-bregman}^{\sym{norm-X}}_{(\Space{V}_{n-1}^{\text{trunc}})^\perp} \left ( \sym{mapping-duality}_{\sym{norm-dualX}} ( \sym{descentdirection}_n ) \right ) ,
			\end{equation}
			where $\Space{V}^\perp$ designates the annihilator of the space $\Space{V}$.
			\end{enumerate}
			\begin{proof}
				For part \ref{itm:solution-manifold-in-intersection} for any $\Vect{z} \in \Space{M}_{\Mat{A}\Vect{x} = \Vect{y}}$ we have to show $\skp{ \sym{searchdirection-orthogonalized}_{n,k} }{ \Vect{z} } = \beta_{n,k}$, which follows directly from the definition of the offsets, \eqref{math:orthogonalization-offsets}.
				Part \ref{itm:dualiterate-in-range-adjoint-operator} follows from
				\begin{align*}
						\sym{mapping-duality}_{\sym{norm-X}} (\sym{iterate}_{n}) - \sym{mapping-duality}_{\sym{norm-X}} (\sym{iterate}_{0}) &= \sym{mapping-duality}_{\sym{norm-X}} (\sym{iterate}_{n}) - \sym{mapping-duality}_{\sym{norm-X}} (\sym{iterate}_{n-1}) + \sym{mapping-duality}_{\sym{norm-X}} (\sym{iterate}_{n-1}) - \ldots + \sym{mapping-duality}_{\sym{norm-X}} (\sym{iterate}_{1}) - \sym{mapping-duality}_{\sym{norm-X}} (\sym{iterate}_0), \\
						&= \sum^{\sym{number-search-directions}_n}_{j=1} \mu_n \sym{searchdirection-orthogonalized}_{n,j} + \ldots + \sum^{\sym{number-search-directions}_1}_{j=1} \mu_1 \sym{searchdirection-orthogonalized}_{1,j},
				\end{align*}
				c.\,f.~\eqref{math:SESOP-update-iterate},
				
				For \ref{itm:searchdirection-projection-intersection} we only need to use the definition~\eqref{math:orthogonalized-search-directions} of $\sym{searchdirection-orthogonalized}_{n, \sym{number-search-directions}_n}$, relations between metric and Bregman projections, see Prop.~\ref{prop:props-bregman-metric-projection}\,\ref{itm:props-bregman-metric-projection-relation} and \cite[Prop.~3.6\,d)]{Schopfer2008}, and equivalencies for Bregman projections, see~\cite[Prop.~3.7\,b)]{Schopfer2008},
				\begin{align*}
					\sym{mapping-duality}_{\sym{norm-dualX}} ( \sym{searchdirection-orthogonalized}_{n,\sym{number-search-directions}_n} ) &= \sym{mapping-duality}_{\sym{norm-dualX}} \left (\sym{descentdirection}_n - \sym{projection-metric}_{\Space{V}_{n-1}^{\text{trunc}}} ( \sym{descentdirection}_n ) \right ) = \sym{mapping-duality}_{\sym{norm-dualX}} \left ( - \sym{projection-bregman}_{\Space{V}_{n-1}^{\text{trunc}} - \sym{descentdirection}_n} (0) \right ) \\
					&= \sym{mapping-duality}_{\sym{norm-dualX}} \left ( \sym{projection-bregman}_{\sym{descentdirection}_n - \Space{V}_{n-1}^{\text{trunc}}} (0) \right ) = \sym{mapping-duality}_{\sym{norm-dualX}} \left ( \sym{projection-bregman}_{\sym{descentdirection}_n + \Space{V}_{n-1}^{\text{trunc}}} (0) \right ) \\
					&= \sym{projection-bregman}_{ \left ( \Space{V}_{n-1}^{\text{trunc}} \right )^\perp} \left ( \sym{mapping-duality}_{\sym{norm-dualX}} ( \sym{descentdirection}_n ) \right ).
				\end{align*}
				
				Part \ref{itm:iterate-projection-intersection} follows as in the proof of \cite[Prop.~4.1]{Schopfer2008}.
			\end{proof}
		\end{lemma}
		
		Before we then may prove weak convergence, we need to show that $\sym{searchdirection-orthogonalized}_{n,\sym{number-search-directions}_n}$ is still a descent direction, where we use the same geometrical arguments as in the generalized Landweber convergence proof, see \cite{Schopfer2006}, on $\sym{descentdirection}_n$ being a descent direction. This is not a straight-forward consequence of this as we have $\Space{U}_n \not\subset \Space{V}_n^{\text{trunc}}$, i.\,e. $\sym{descentdirection}_n$ is not generally contained in $\Space{V}_n^{\text{trunc}}$, c.\,f.~\eqref{math:orthogonalized-search-directions-projections}.
		\pagebreak[0]\begin{lemma}[Descent direction property]\label{lem:descent-direction-property}
			Any $\sym{searchdirection-orthogonalized}_{n,\sym{number-search-directions}_n}$ resulting from \eqref{math:orthogonalized-search-directions} is always a descent direction, i.\,e.~there is a $\sym{stepwidth}_n \in \reel^{\sym{number-search-directions}_n}$ and $S_n > 0$ with
			\begin{equation}
				h_n(\sym{stepwidth}_n) \leq h_n(0) - S_n,
			\end{equation}
			which for any $\Vect{z} \in \Space{M}_{\Mat{A}\Vect{x} = \Vect{y}}$ can also be written as
			\begin{equation}
				\label{math:descent-direction-property}
				\sym{distance-bregman}_{\Space{X}} (\sym{iterate}_{n+1}, \Vect{z}) \leq \sym{distance-bregman}_{\Space{X}} (\sym{iterate}_{n}, \Vect{z}) - S_n.
			\end{equation}
			\begin{proof}
				We assume $\sym{residual}_n \neq 0$ and $\sym{iterate}_n \neq 0$.
				We set
				\begin{equation}
					\widetilde{\sym{stepwidth}}_n := (0,\ldots, 0, \nu_n) \quad \text{with } \nu_n := \frac{ \tau_n \norm{ \sym{iterate}_n }_{\Space{X}}^{p-1} }{ \norm{ \sym{searchdirection-orthogonalized}_{n,\sym{number-search-directions}_n} }_{\DualSpace{X}} } 
				\end{equation}
				with $\tau_n \in (0,1]$ chosen as to fulfil with $\gamma \in (0,1)$, c.\,f.~Theorem~\ref{thm:xu-roach-inequality-smoothness},
				\begin{equation*}
					\frac{ \sym{modulus-smoothness}_{\DualSpace{X}} ( \tau_n ) }{ \tau_n } = \sym{modulus-smoothness}_{\DualSpace{X}} ( 1 ) \wedge \left ( \frac{ \gamma }{ 2^{\sym{norm-dualX}} G_{\sym{norm-X}} } \frac{ \sym{residual}_n^r }{ \norm{ \sym{iterate}_n }_{\Space{X}} \norm{ \sym{searchdirection-orthogonalized}_n }_{\DualSpace{X}} } \right ).
				\end{equation*}
				Let $\sym{stepwidth}_n = \argmin_{t \in \reel^n} h_n(t)$, then we estimate with the Xu-Roach inequality~\eqref{math:xu-roach-inequality-smoothness} and using $\widetilde{\sym{stepwidth}}_n$,
				\begin{align*}
					h_n(\sym{stepwidth}_n) &\leq h_n (\widetilde{\sym{stepwidth}}_n) \\
					&\leq \tfrac 1 {\sym{norm-dualX}} \norm{ \sym{iterate}_n }_{\Space{X}}^{\sym{norm-X}} - \nu_n \skp{ \sym{searchdirection-orthogonalized}_{n,\sym{number-search-directions}_n} }{ \sym{iterate}_n } + \frac 1 {\sym{norm-dualX}} \widetilde{\sigma} \left ( \sym{mapping-duality} ( \sym{iterate}_n), \nu_n \sym{searchdirection-orthogonalized}_{n,\sym{number-search-directions}_n} \right ) + \skp{ \sym{searchdirection-orthogonalized-precursor}_{n,\sym{number-search-directions}_n} }{ \Vect{y} } \\
					&= \tfrac 1 {\sym{norm-dualX}} \norm{ \sym{iterate}_n }_{\Space{X}}^{\sym{norm-X}} - \nu_n \skp{ \sym{descentdirection-precursor}_n - \sum^{\sym{number-search-directions}_{n-1}}_{k=1} \sym{orthogonalization-coefficient}_{n,k} \sym{searchdirection-orthogonalized-precursor}_{n-1,k} }{ \sym{residualdirection}_n }  + \frac 1 {\sym{norm-dualX}} \widetilde{\sigma} \left ( \sym{mapping-duality} ( \sym{iterate}_n), \nu_n \sym{searchdirection-orthogonalized}_{n,\sym{number-search-directions}_n} \right ) \\
					&= \tfrac 1 {\sym{norm-dualX}} \norm{ \sym{iterate}_n }_{\Space{X}}^{\sym{norm-X}} - \nu_n \sym{residual}_n^r + \frac 1 {\sym{norm-dualX}} \widetilde{\sigma} \left ( \sym{mapping-duality} ( \sym{iterate}_n), \nu_n \sym{searchdirection-orthogonalized}_{n,\sym{number-search-directions}_n} \right ) \\
					&= \tfrac 1 {\sym{norm-dualX}} \norm{ \sym{iterate}_n }_{\Space{X}}^{\sym{norm-X}} - \tau_n \norm{ \sym{iterate}_n }^{p-1} \frac{ \sym{residual}_n^r }{ \norm{ \sym{searchdirection-orthogonalized}_{n,\sym{number-search-directions}_n} }_{\DualSpace{X}} } + \frac 1 {\sym{norm-dualX}} \widetilde{\sigma} \left ( \sym{mapping-duality} ( \sym{iterate}_n), \nu_n \sym{searchdirection-orthogonalized}_{n,\sym{number-search-directions}_n} \right ),
				\end{align*}
				where we have used the orthogonality stated in Corollary~\ref{cor:expanding-subspace}\,\ref{item:expanding-subspace-orthogonal-residual}.
				As the metric projection is non-expanding, c.\,f.~\eqref{math:orthogonalized-search-directions-projections}, we have
				\begin{equation*}
					\norm{ \sym{searchdirection-orthogonalized}_{n,\sym{number-search-directions}_n} }_{\DualSpace{X}} \leq \norm{ \dual{\Mat{A}} \sym{mapping-duality}_{\sym{norm-Y}} (\Mat{A} \sym{iterate}_n - \Vect{y} )} \leq \norm{ \Mat{A} } \sym{residual}_n^{r-1}.
				\end{equation*}
				Lemma~\ref{lem:sigma-upperbound} allows us to bound the last summand using the requirement on $\tau_n$,
				\begin{align*}
					\tfrac 1 {\sym{norm-dualX}} \widetilde{\sigma}_{\sym{norm-dualX}} \left ( \sym{mapping-duality} ( \sym{iterate}_n), \nu_n \sym{searchdirection-orthogonalized}_{n,\sym{number-search-directions}_n} \right ) &\leq 2^{\sym{norm-dualX}} G_{\sym{norm-X}} \norm{ \sym{iterate}_n }_{\Space{X}}^{\sym{norm-X}} \sym{modulus-smoothness}_{\DualSpace{X}} (\tau_n) \\
					&\leq \tau_n 2^{\sym{norm-dualX}} G_{\sym{norm-X}} \norm{ \sym{iterate}_n }_{\Space{X}}^{\sym{norm-X}} \frac{ \sym{modulus-smoothness}_{\DualSpace{X}} (\tau_n) }{\tau_n} \\
					&\leq \tau_n 2^{\sym{norm-dualX}} G_{\sym{norm-X}} \norm{ \sym{iterate}_n }_{\Space{X}} ^{\sym{norm-X}} \frac{ \gamma }{ 2^{\sym{norm-dualX}} G_{\sym{norm-X}} } \frac{ \sym{residual}_n^r }{ \norm{ \sym{iterate}_n }_{\Space{X}} \norm{ \sym{searchdirection-orthogonalized}_{n,\sym{number-search-directions}_n} }_{\DualSpace{X}} } \\
					&= \gamma \tau_n \norm{ \sym{iterate}_n }_{\Space{X}}^{p-1}  \frac{ \sym{residual}_n^r }{ \norm{ \sym{searchdirection-orthogonalized}_{n,\sym{number-search-directions}_n} }_{\DualSpace{X}} }.
				\end{align*}
				All together we obtain
				\begin{align}
					h_n(\sym{stepwidth}_n) &\leq \tfrac 1 {\sym{norm-dualX}} \norm{ \sym{iterate}_n }_{\Space{X}}^{\sym{norm-X}} - (1 - \gamma) \tau_n \norm{ \sym{iterate}_n }_{\Space{X}}^{p-1} \frac{ \sym{residual}_n^r }{ \norm{ \sym{searchdirection-orthogonalized}_{n,\sym{number-search-directions}_n} } } \nonumber \\
					\label{math:bound-bregman-distance-residual}
					&\leq \tfrac 1 {\sym{norm-dualX}} \norm{ \sym{iterate}_n }_{\Space{X}}^{\sym{norm-X}} - \frac{ (1 - \gamma) }{ \norm{\Mat{A} } } \tau_n \norm{ \sym{iterate}_n }_{\Space{X}}^{p-1} \sym{residual}_n.
				\end{align}
				We are done as all factors in $S_n := \tfrac{ (1 - \gamma) }{ \norm{\Mat{A} } } \tau_n \norm{ \sym{iterate}_n }_{\Space{X}}^{p-1} \sym{residual}_n$ are positive and knowing $h_n(0) 
				= \tfrac 1 {\sym{norm-dualX}} \norm{ \sym{iterate}_n }_{\Space{X}}^{\sym{norm-X}}$.
				
				Lastly, let us show how $h_n(\sym{stepwidth})$ and $\sym{distance-bregman}_{\sym{norm-X}} (\sym{iterate}_{n+1}(\sym{stepwidth}), \Vect{z})$ are related using \eqref{math:distance-bregman-identity-power} and \eqref{math:SESOP-update-iterate}. For any $\Vect{z} \in \Space{M}_{\Mat{A}\Vect{x} = \Vect{y}}$ we have
				\begin{align*}
					\sym{distance-bregman}_{\sym{norm-X}} (\sym{iterate}_{n+1}(\sym{stepwidth}), \Vect{z}) &= \tfrac 1 {\sym{norm-dualX}} \norm{ \sym{iterate}_{n+1}(\sym{stepwidth}) }_{\Space{X}}^{\sym{norm-X}} - \skp{ \sym{mapping-duality}_{\sym{norm-X}} (\sym{iterate}_n) - \sum_{k=1}^{\sym{number-search-directions}_n} \sym{stepwidth}_k \sym{searchdirection-orthogonalized}_{n,k} }{ \Vect{z} } + \tfrac 1 p \norm{ \Vect{z} }_{\Space{X}}^{\sym{norm-X}} \\
					&= \tfrac 1 {\sym{norm-dualX}} \norm{ \sym{iterate}_{n+1}(\sym{stepwidth}) }_{\Space{X}}^{\sym{norm-X}} - \skp{ \sym{mapping-duality}_{\sym{norm-X}} (\sym{iterate}_n) }{ \Vect{z} } + \sum_{k=1}^{\sym{number-search-directions}_n} \sym{stepwidth}_k \beta_{n,k} + \tfrac 1 p \norm{ \Vect{z} }_{\Space{X}}^{\sym{norm-X}} \\
					&= h_n(\sym{stepwidth}) - \skp{ \sym{mapping-duality}_{\sym{norm-X}}  (\sym{iterate}_n) }{ \Vect{z} } + \tfrac 1 p \norm{ \Vect{z} }_{\Space{X}}^{\sym{norm-X}}.
				\end{align*}
				As the last two terms are constant with respect to $\sym{stepwidth}$, they cancel out when considering the difference $h_n(0) - h_n(\sym{stepwidth})$. Hence, \eqref{math:descent-direction-property} holds.
			\end{proof}
		\end{lemma}
		
		Finally, we arrive at the proof of convergence. 


		\pagebreak[0]\begin{theorem}[Weak convergence, truncated, semi-orthogonalized search space]\label{thm:convergence-truncated-orthogonalized-searchspace}
			Given a uniformly convex and smooth Banach space $\Space{X}$ with sequentially weak-to-weak continuous duality mapping and an arbitrary Banach space $\Space{Y}$, 
			then with $1 \leq \sym{number-search-directions}_n \leq n$ and search space $\Space{V}_n^{\text{trunc}}$, given by \eqref{math:orthogonalized-search-directions-olddirections} and \eqref{math:orthogonalized-search-directions}, 
			method~\ref{method:sesop} either stops after a finite number $n \in \nat$ of iterations (in case $\sym{residual}_n=0$) with $\sym{iterate}_n$ being the Bregman projection $\sym{truesolution} = \sym{projection-bregman}_{\Space{M}_{\Mat{A}\Vect{x}=\Vect{y}}} (\sym{iterate}_0)$ of $\sym{iterate}_0$ onto the solution manifold $\Space{M}_{\Mat{A}\Vect{x}=\Vect{y}}$ or the sequence of the iterates $\{\sym{iterate}_n\}_n$ converges weakly to $\sym{truesolution}$.
			\begin{proof}
				In case $\sym{residual}_{n^\ast} = 0$ for some $n^\ast$, then we have $\sym{iterate}_{n^\ast} \in \Space{M}_{\Mat{A}\Vect{x}=\Vect{y}}$ and we are done by \cite[Proposition~3.7\,b)]{Schopfer2008} together with the optimality condition in Lemma~\ref{lem:condition-optimality}\,\ref{itm:minimum-norm-solution-optimality}.
				
				Let us then assume $\sym{residual}_{n} \neq 0$ for all $n$.
				Lemma~\ref{lem:descent-direction-property} ensures that $\{\sym{distance-bregman}_{\sym{norm-X}} (\sym{iterate}_n, \Vect{z}) \}_n$ for all $\Vect{z} \in \Space{M}_{\Mat{A}\Vect{x}=\Vect{y}}$ is strictly decreasing. Therefore, $\{\sym{distance-bregman}_{\sym{norm-X}} (\sym{iterate}_n, \Vect{z}) \}_n$ is bounded from above by $\{\sym{distance-bregman}_{\sym{norm-X}} (\sym{iterate}_0, \Vect{z}) \}_n$.
				
				Then, Proposition~\ref{prop:distance-bregman-properties} \ref{itm:bregman-distance-boundedness} assures that $\{ \sym{iterate}_n \}_n$ is bounded. As we require $\Space{X}$ to be uniformly convex and hence reflexive by the Milman-Pettis theorem, \cite[\citesection~II.2, Thm.~2.9]{Cioranescu1990}, every subsequence of $\{ \sym{iterate}_n \}_n$ has in turn a subsequence $\{ \sym{iterate}_{n_k} \}_k$ that converges weakly to some $\widetilde{\Vect{x}} \in \Space{X}$, see \cite[\citechapter~8,Thm~4.2]{Schechter1971}. 

				The proof of $\{ \sym{residual}_{n_k} \}_k$ being a null sequence follows in exactly the same way as in \cite[\citepage~320]{Schopfer2006}.
				Then, we even have $\widetilde{\Vect{x}} \in \Space{M}_{\Mat{A}\Vect{x}=\Vect{y}}$ as $\sym{residual}_{n_k} = \norm{\Mat{A}\sym{iterate}_{n_k} - \Vect{y}}_{\Space{Y}} \rightarrow 0$ with $k \rightarrow \infty$, i.\,e.~$\Mat{A}\widetilde{\Vect{x}} = \Vect{y}$.
				
				As $\overline{\Range{\dual{\Mat{A}}}}$ is convex and norm-closed, it is also weakly closed, see \cite[\citechapter~5, Thm.~3.13]{Dunford1957}. This together with Lemma~\ref{lem:intersection-hyperplanes} \ref{itm:dualiterate-in-range-adjoint-operator} implies $\sym{mapping-duality}_{\sym{norm-X}} (\widetilde{\Vect{x}}) - \sym{mapping-duality}_{\sym{norm-X}} (\sym{iterate}_0) \in \overline{\Range{\dual{\Mat{A}}}}$ given $\sym{mapping-duality}_{\sym{norm-X}}$ is sequentially weak-to-weak continuous, see the note below. By the optimality condition Lemma~\ref{lem:condition-optimality}\,\ref{itm:minimum-norm-solution-optimality} with the requirement $\sym{mapping-duality}_{\sym{norm-X}}(\sym{iterate}_0) \in \overline{\Range{\dual{\Mat{A}}}}$ and with $\Vect{z} + \Null{\Mat{A}} = \Space{M}_{\Mat{A}\Vect{x}=\Vect{y}}$ for all $\Vect{z} \in \Space{M}_{\Mat{A}\Vect{x}=\Vect{y}}$, we then conclude $\widetilde{\Vect{x}} = \sym{projection-bregman}_{\Space{M}_{\Mat{A}\Vect{x}=\Vect{y}}} (\sym{iterate}_0)$.
				
				As we have shown that for every subsequence of $\{\sym{iterate}_n \}_n$ there is a subsequence that in turn converges weakly to the same limit $\sym{projection-bregman}_{\Space{M}_{\Mat{A}\Vect{x}=\Vect{y}}} (\sym{iterate}_0)$, then this is valid for the sequence $\{\sym{iterate}_n \}_n$, too, see \cite[\citesection~10.5]{Zeidler1986}.
			\end{proof}
		\end{theorem}
		
		Note that the duality mappings of $\ell_{\sym{norm-X}}$-spaces, where $1<p<\infty$, are sequentially weak-to-weak-contin\-u\-ous, see \cite[Remark~4.3]{Schopfer2008}.
				
\section{Experiments}\label{sec:experiments}

	In the following we will perform numerical experiments that compare method~\ref{method:sesop} with search spaces $\Space{U}^{\text{trunc}}_n$ and $\Space{V}^{\text{trunc}}_n$. The first part is identical to the experiments in \cite[\citesection~5]{Schopfer2008}, namely solving $\Mat{A} \Vect{x} = \Vect{y}$ for various $\ell_{\sym{norm-X}}$ spaces and uniformly distributed random matrices $\Mat{A}$ and right-hand sides $\Vect{y}$. In the second part we solve inverse problems in 2D computerized tomography with the Radon transform as linear operator.
	
	All experiments habe been performed on a single core of an Intel Xeon E3-1270 cpu with 3.50GHz.
	Note that \glostext{SESOP} and the above described orthogonalization procedure have been implemented in the C++ library \emph{BASSO} (\emph{BA}nach \emph{S}equential \emph{S}ubspace \emph{O}ptimizer), based on the Eigen3 library~\cite{eigen3} for the linear algebra routines. It is available on request from the authors.

	\subsection{Toy Problem}\label{sec:experiments-toyproblem}
	
		We first look at the inverse toy problem of a random matrix and a random right-hand side to be formally inverted as is done frequently in the literature, but with well-known short-comings, see \cite{Loris2009}.

		Most importantly in this experiment, we want to check the case of $p=2$ for a single and for multiple search directions. With \glostext{SESOP} there is still a significant decrease in required iterations from single to multiple search direction for the $\ell_2$ space. With \glostext{CG} no such difference arises due to the conjugacy property. 
		Furthermore, to assess a possible speed-up of the orthogonalized search directions, we look at various $\ell_{\sym{norm-X}}$ spaces and norms.
	
		\subsubsection{Procedure}\label{sec:experiments-toyproblem-procedure}
	
			To this end, we create a uniformly distributed random matrix $\sym{matrix-forward} \in [-1,1]^{\sym{dimension-rhs} \times \sym{dimension-solution}}$ with $\sym{dimension-rhs}=1000$ and $\sym{dimension-solution}=5000$, representing a discretized version of some random operator.
		
			Next, we want to create a "solution" $\sym{truesolution} \in \reel^{\sym{dimension-solution}}$ to a random right-hand side $y \in \reel^{\sym{dimension-rhs}}$ in the sense that this solution should be a minimum-norm solution in an $\ell_{\sym{norm-X}}$-space with $\sym{norm-X} \in \{1.1, 1.2, 1.5, 2, 3, 6, 10\}$. Therefore, we create a random right-hand side precursor $y^\ast \in [-1,1]^{\sym{dimension-rhs}}$ and calculate the minimum-norm solution $x^\dag$ as follows
			\begin{equation}\label{math:minimum-norm-solution-rhs}
				\sym{truesolution} := \frac {\sym{mapping-duality}^\adj_{\sym{norm-dualX}}( \dual{\Mat{A}} y^\ast) } { \norm{\sym{mapping-duality}^\adj_{\sym{norm-dualX}}( \dual{\Mat{A}} y^\ast)}_{\sym{norm-X}} }
			\end{equation}
			And 	then we finally obtain the right-hand side as $\reel^{\sym{dimension-rhs}} \ni y = \sym{matrix-forward} \sym{truesolution}$.
			
			We use 10 different seeds $s \in \{420, \ldots, 429\}$ for the random number generator and calculate average iteration counts $n$ and standard deviations $\sigma_n$ over these 10 runs with otherwise identical parameters. The iteration is stopped at either $n > 20.000$ or if the relative residual $\frac{ \norm{ \Mat{A}\Vect{x}_n - \Vect{y} }_{\Space{Y}} }{ \norm{ \Vect{y} }_{\Space{Y}} }$ drops below $10^{-4}$.
			
			Note that we always set $\Space{Y}=\ell_2(\mathbb{R}^{\sym{dimension-rhs}})$.
	
		\subsubsection{SESOP}\label{sec:experiments-toyproblem-SESOP}
		
			First, we reproduce the results from \cite[\citesection~5]{Schopfer2008} to elucidate any possible differences that arise from different implementations, see Table~\ref{tab:solution_iterations-reproducingresults}. There, the matrix dimension is $1000 \times 5000$, the $\ell_{\sym{norm-X}}$-norms of $\Space{X}$ are given by $\sym{norm-X} \in \{1.2, 1.5, 6, 10\}$, the power of the gauge function of the duality mapping $\sym{mapping-duality}_{\sym{norm-X}}$ is given by $\begin{cases}\sym{norm-X}, \sym{norm-X} \geq 2 \\ 2, \text{else}\end{cases}$ and we use search direction numbers $\sym{number-search-directions} \in \{2,4,6\}$. 
	
			\begin{table}
				\begin{center}
				\def\mydim{5000}%
				\caption{Comparison of \protect\glostext{SESOP} implementations from \cite{Schopfer2008} and this work (without orthogonalization): Average iteration counts $n$ and standard deviations $\sigma_n$ for various $\protect\sym{norm-X}$ and $\protect\sym{number-search-directions}$ values and matrix dimension $m=\mydim$.}
				\label{tab:solution_iterations-reproducingresults}
				\subtable[from \cite{Schopfer2008}\label{tab:solution_iterations-Schoepfer2008}]{
					\input{./figures/Schoepfer2008-FunctionFit-results}
				} \goodgap
				\subtable[This work\label{tab:solution_iterations-thiswork}]{
					\def\mysuffix{_unprojected_small}%
					\input{./figures/solution_iterations_small}
				}
				\end{center}
			\end{table}

			Taking into account that in Table~\ref{tab:solution_iterations-Schoepfer2008} only a single run is given and comparing this to our averages and standard deviations, the iteration counts are in very good agreement up until $\sym{norm-X}=2$. For larger $\sym{norm-X}$ the discrepancy is quite large. However, this is also true for the standard deviations.
			
			Hence, overall we do not see any discrepancy resulting from the different implementations, i.\,e.~we have a solid base for comparing the results with the original MatLab implementation of \cite{Schopfer2008}.
			
		\subsubsection{Orthogonalized Search Directions}\label{sec:experiments-toyproblem-orthogonalized-directions}
		
			Next, we look at the change in iteration counts between the search space $\Space{U}^{\text{trunc}}_n$ used in \glostext{SESOP} and the search spaces $\Space{V}_n^{\text{trunc}}$ using metric projections proposed in this article. We use various numbers of search directions $\sym{number-search-directions} \in \{ 1,2,4,6\}$ and a full exemplary range of $\ell_{\sym{norm-X}}$ spaces with $\sym{norm-X} \in \{1.1, 1.2, 1.5, 2, 3, 6, 10\}$ and the power type of the duality mapping chosen as before. 

			\begin{table}
				\def\mydim{5000}%
				\caption{Averaged iteration counts $n$ and standard deviations $\sigma_n$ for various $p$ and $N$ values and $m=\mydim$ for solving with \protect\glostext{SESOP} for the minimum-norm solution of $\Mat{A}\Vect{x}=\Vect{y}$: \emph{unorth}ogonalized refers to $\Space{U}^{\text{trunc}}$ and \emph{metric} projection to $\Space{V}^{\text{trunc}}$.}
				\label{tab:solution_iterations-comparingprojected}
				\begin{center}
					\ifdefined\mydim\empty\else\def\mydim{5000}\fi%

\pgfplotstableread[
	columns={p,iterations},
	col sep=comma,
]{./figures/pnorm-iterations-Competition_of_search_spaces-dim_\mydim-orthogonal_directions_1.csv}\metriciterationtable%
\pgfplotstableread[
	columns={p,iterations},
	col sep=comma,
]{./figures/pnorm-iterations-Competition_of_search_spaces-dim_\mydim-orthogonal_directions_0.csv}\iterationtable%

\pgfplotstableset{
	create on use/p/.style={create col/copy column from table={\iterationtable}{p}},
	create on use/N/.style={create col/copy column from table={\iterationtable}{N}},
	create on use/avg_iterations_metric/.style={create col/copy column from table={\metriciterationtable}{avg_iterations}},
	create on use/var_iterations_metric/.style={create col/copy column from table={\metriciterationtable}{var_iterations}},
	create on use/avg_iterations/.style={create col/copy column from table={\iterationtable}{avg_iterations}},
	create on use/var_iterations/.style={create col/copy column from table={\iterationtable}{var_iterations}},
	create on use/err_iterations_metric/.style={
		create col/expr={sqrt(\thisrow{var_iterations_metric}-\thisrow{avg_iterations_metric}*\thisrow{avg_iterations_metric})}
	},
	create on use/err_iterations/.style={
		create col/expr={sqrt(\thisrow{var_iterations}-\thisrow{avg_iterations}*\thisrow{avg_iterations})}
	},
}

\pgfplotstablegetrowsof\iterationtable
\edef\myrows{\pgfplotsretval}
\pgfplotstablenew[
	columns={p,N,avg_iterations_metric,var_iterations_metric,avg_iterations,var_iterations},
]{\myrows}\combinedtable

\pgfplotstabletypeset[
	columns={p,N,avg_iterations,err_iterations,avg_iterations_metric,err_iterations_metric},
	col sep=&,
	row sep=\\,
	every first row/.style={before row={\rowcolor[gray]{0.9}}},
	every nth row={8[0]}{before row={\rowcolor[gray]{0.9}}},
	every nth row={8[+1]}{before row={\rowcolor[gray]{0.9}}},
	every nth row={8[+2]}{before row={\rowcolor[gray]{0.9}}},
	every nth row={8[+3]}{before row={\rowcolor[gray]{0.9}}},
	every head row/.style={	
		before row=\hline,
		after row=\hline\hline,
	},
	every last row/.style={after row=\hline},
	columns/p/.style={
		column type/.add={|}{}
	},
	columns/N/.style={
		column type/.add={}{|}
	},
	columns/avg_iterations/.style={
		column name={$n^{\text{unorth}}$},
		sci subscript,
	},
	columns/err_iterations/.style={
		column name={$\sigma^{\text{unorth}}_n$},
		column type/.add={}{|},
		sci subscript,
	},
	columns/avg_iterations_metric/.style={
		column name={$n^{\text{metric}}$},
		sci subscript,
	},
	columns/err_iterations_metric/.style={
		column name={$\sigma^{\text{metric}}_n$},
		column type/.add={}{|},
		sci subscript,
	},
	every row 0 column p/.style={column type/.add={|}{|}},
	every row 0 column err_iterations/.style={column type/.add={}{|}},
]\combinedtable

				\end{center}
			\end{table}
	
			The average iteration counts and standard deviations are given in Table~\ref{tab:solution_iterations-comparingprojected}.
			
			First of all, we notice that with $\sym{norm-X} \in \{1.2,1.5,2\}$ there is not much change in the average iteration count $n$ over different number of search directions $\sym{number-search-directions}$. This holds on for $\sym{norm-X}=1.5$ within output precision of $10^{-7}$ and for $\sym{norm-X}=2$ within full floating point numerical precision. Hence, we see that when $\Space{X} = \ell_2 \left (\reel^{\sym{dimension-solution}} \right )$ is a Hilbert-space, the orthogonality is maintained between all search directions from $\sym{searchdirection-orthogonalized}_{0,\sym{number-search-directions}_0}$ up to $\sym{searchdirection-orthogonalized}_{n,\sym{number-search-directions}_n}$. This is the behavior expected from \glostext{CG} methods and we elucidate this further in the next section. It is maintained to some extent when $\sym{norm-X}$ is close to $2$.

			Second, for $\sym{norm-X} \in \{1.5, 2\}$ we observe that for large number of search directions $\sym{number-search-directions}$ the average number of iteration steps $n$ becomes equivalent for both of the search spaces $\Space{U}^{\text{trunc}}_n$ and $\Space{V}_n^{\text{trunc}}$. This indicates that using more than one search direction, the central idea of \cite{Schopfer2008}, is indeed substantial.
			
			And last, iteration counts and deviations become very large for $\sym{norm-X} \rightarrow 1$ and $\sym{norm-X} \rightarrow \infty$, i.\,e.~when the $\ell_{\sym{norm-X}}$ spaces are no longer smooth.
	
			\begin{figure}[htbp]
				\def\myN{1}%
				\subfigure[Iterations\label{fig:solution_iterations-comparingprojected_iterations}]{
					\resizebox{0.45\textwidth}{!}{
						\ifdefined\Space\empty\else\newcommand{\Space}[1]{#1}\fi%

\ifdefined\mylowerylimit\empty\else\def\mylowerylimit{6e0}\fi%
\ifdefined\myupperylimit\empty\else\def\myupperylimit{5e3}\fi%
\ifdefined\myN\empty\else\def\myN{1}\fi%
\ifdefined\mydim\empty\else\def\mydim{5000}\fi%
\ifdefined\mypath\empty\else\newcommand{\mypath}{/home/heber/Documents/Experiments/BregmanProjectedSearchDirections}\fi%

\pgfplotsset{
    discard if not/.style 2 args={
        x filter/.code={%
            \edef\tempa{\thisrow{#1}}%
            \edef\tempb{#2}%
            \ifx\tempa\tempb%
            \else%
                \def\pgfmathresult{inf}%
            \fi%
        }
    }
}

\begin{tikzpicture}[]
\begin{loglogaxis}[
	ymin=\mylowerylimit,
	ymax=\myupperylimit,
	grid=major,
	xtick={1.2,1.5,2,3,6,10},
	extra x ticks={1.1},
	extra x tick labels={1.1},
	major tick length={10pt},
	extra x tick style={rotate=0,major tick length=20pt,xtick pos=left, xtick align=center},
	ymajorgrids=true,
	log ticks with fixed point,
	xticklabel shift=2pt,
	legend pos=south east,
	legend columns=1,
	ylabel={iterations $n$},
	xlabel={$\ell_p$-norm}
]
\addplot+[
 	thick,
	discard if not={N}{\myN},
	error bars/.cd,y dir=both,y explicit,
] table[x=p,y=avg_iterations,y error expr={sqrt(\thisrow{var_iterations}-\thisrow{avg_iterations}^2)},col sep=comma]{\mypath /RandomMatrix/data/pnorm-iterations-Competition_of_search_spaces-dim_\mydim-orthogonal_directions_0.csv};
\addlegendentry{$\Space{U}^{\text{trunc}}_n$}

\addplot+[
 	thick,
	discard if not={N}{\myN},
	error bars/.cd,y dir=both,y explicit,
] table[x=p,y=avg_iterations,y error expr={sqrt(\thisrow{var_iterations}-\thisrow{avg_iterations}^2)},col sep=comma,]{\mypath /RandomMatrix/data/pnorm-iterations-Competition_of_search_spaces-dim_\mydim-orthogonal_directions_1.csv};
\addlegendentry{$\Space{V}^{\text{trunc}}_n$}

\ifdefined\mytitle{%
\node[rectangle,draw,thick,fill=white] (title) at (1,5) {\mytitle};
}\fi%

\end{loglogaxis}
\end{tikzpicture}
					}
				}
				\subfigure[Runtimes\label{fig:solution_iterations-comparingprojected_runtimes}]{
					\resizebox{0.45\textwidth}{!}{
						\ifdefined\Space\empty\else\newcommand{\Space}[1]{#1}\fi%

\ifdefined\mylowerylimit\empty\else\def\mylowerylimit{1e-1}\fi%
\ifdefined\myupperylimit\empty\else\def\myupperylimit{5e2}\fi%
\ifdefined\myN\empty\else\def\myN{1}\fi%
\ifdefined\mydim\empty\else\def\mydim{5000}\fi%
\ifdefined\mypath\empty\else\newcommand{\mypath}{/home/heber/Documents/Experiments/BregmanProjectedSearchDirections}\fi%

\pgfplotsset{
    discard if not/.style 2 args={
        x filter/.code={%
            \edef\tempa{\thisrow{#1}}%
            \edef\tempb{#2}%
            \ifx\tempa\tempb%
            \else%
                \def\pgfmathresult{inf}%
            \fi%
        }
    }
}

\begin{tikzpicture}[]
\begin{loglogaxis}[
	ymin=\mylowerylimit,
	ymax=\myupperylimit,
	grid=major,
	xtick={1.2,1.5,2,3,6,10},
	extra x ticks={1.1},
	extra x tick labels={1.1},
	major tick length={10pt},
	extra x tick style={rotate=0,major tick length=20pt,xtick pos=left, xtick align=center},
	ymajorgrids=true,
	log ticks with fixed point,
	xticklabel shift=2pt,
	legend pos=south east,
	legend columns=1,
	ylabel={runtime [s]},
	xlabel={$\ell_p$-norm}
]
\addplot+[
 	thick,
	discard if not={N}{\myN},
	error bars/.cd,y dir=both,y explicit,
] table[x=p,y=avg_runtime,y error expr={sqrt(\thisrow{var_runtime}-\thisrow{avg_runtime}^2)},col sep=comma]{\mypath /RandomMatrix/data/pnorm-runtimes-Competition_of_search_spaces-dim_\mydim-orthogonal_directions_0.csv};
\addlegendentry{$\Space{U}^{\text{trunc}}_n$}

\addplot+[
 	thick,
	discard if not={N}{\myN},
	error bars/.cd,y dir=both,y explicit,
] table[x=p,y=avg_runtime,y error expr={sqrt(\thisrow{var_runtime}-\thisrow{avg_runtime}^2)},col sep=comma,]{\mypath /RandomMatrix/data/pnorm-runtimes-Competition_of_search_spaces-dim_\mydim-orthogonal_directions_1.csv};
\addlegendentry{$\Space{V}^{\text{trunc}}_n$}

\ifdefined\mytitle{%
\node[rectangle,draw,thick,fill=white] (title) at (1,5) {\mytitle};
}\fi%

\end{loglogaxis}
\end{tikzpicture}
					}
				}
				\caption{Iterations and runtimes solving with \protect\glostext{SESOP} and two different search spaces for the minimum-norm solution of $\Mat{A}\Vect{x}=\Vect{y}$ with uniformly random $\Mat{A} \in \reel^{1000 \times 5000}$ and a single search direction.}
				\label{fig:solution_iterations-comparingprojected}
			\end{figure}
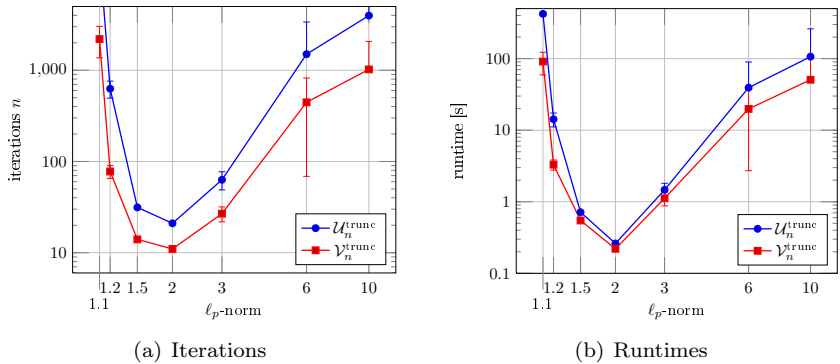

			In Figure~\ref{fig:solution_iterations-comparingprojected} we depict for a single search direction both iteration counts and the total runtime for solving for the minimum-norm solution up to a relative residual threshold of $10^{-4}$ or up to $20.000$ iteration steps. Here, we want to compare the method's performance with either search space directly.
			
			We notice that iteration counts for the orthogonalized search space $\Space{V}^{\text{trunc}}_n$ are at least a factor of 2-3 below the ones for the search space $\Space{U}^{\text{trunc}}_n$. This holds over all values of $\sym{norm-X}$.
			This reduced number of iterations required for the same residual threshold is the reason, why despite the additional computational effort for the orthogonalization, also the runtimes show the same trend between the compared search spaces up to a similar factor.
			

			\begin{figure}[htbp]
			\def\mysuffix{_correctprojected}%
				\resizebox{.45\textwidth}{!}{
					\def\myN{1}%
					\edef\mytitle{$N=\myN$}%
					\ifdefined\Space\empty\else\newcommand{\Space}[1]{#1}\fi%

\ifdefined\mylowerylimit\empty\else\def\mylowerylimit{1e-1}\fi%
\ifdefined\myupperylimit\empty\else\def\myupperylimit{5e2}\fi%
\ifdefined\myN\empty\else\def\myN{1}\fi%
\ifdefined\mydim\empty\else\def\mydim{5000}\fi%
\ifdefined\mypath\empty\else\newcommand{\mypath}{/home/heber/Documents/Experiments/BregmanProjectedSearchDirections}\fi%

\pgfplotsset{
    discard if not/.style 2 args={
        x filter/.code={%
            \edef\tempa{\thisrow{#1}}%
            \edef\tempb{#2}%
            \ifx\tempa\tempb%
            \else%
                \def\pgfmathresult{inf}%
            \fi%
        }
    }
}

\begin{tikzpicture}[]
\begin{loglogaxis}[
	ymin=\mylowerylimit,
	ymax=\myupperylimit,
	grid=major,
	xtick={1.2,1.5,2,3,6,10},
	extra x ticks={1.1},
	extra x tick labels={1.1},
	major tick length={10pt},
	extra x tick style={rotate=0,major tick length=20pt,xtick pos=left, xtick align=center},
	ymajorgrids=true,
	log ticks with fixed point,
	xticklabel shift=2pt,
	legend pos=south east,
	legend columns=1,
	ylabel={runtime [s]},
	xlabel={$\ell_p$-norm}
]
\addplot+[
 	thick,
	discard if not={N}{\myN},
	error bars/.cd,y dir=both,y explicit,
] table[x=p,y=avg_runtime,y error expr={sqrt(\thisrow{var_runtime}-\thisrow{avg_runtime}^2)},col sep=comma]{\mypath /RandomMatrix/data/pnorm-runtimes-Competition_of_search_spaces-dim_\mydim-orthogonal_directions_0.csv};
\addlegendentry{$\Space{U}^{\text{trunc}}_n$}

\addplot+[
 	thick,
	discard if not={N}{\myN},
	error bars/.cd,y dir=both,y explicit,
] table[x=p,y=avg_runtime,y error expr={sqrt(\thisrow{var_runtime}-\thisrow{avg_runtime}^2)},col sep=comma,]{\mypath /RandomMatrix/data/pnorm-runtimes-Competition_of_search_spaces-dim_\mydim-orthogonal_directions_1.csv};
\addlegendentry{$\Space{V}^{\text{trunc}}_n$}

\ifdefined\mytitle{%
\node[rectangle,draw,thick,fill=white] (title) at (1,5) {\mytitle};
}\fi%

\end{loglogaxis}
\end{tikzpicture}
				}
				\resizebox{.45\textwidth}{!}{
					\def\myN{4}%
					\edef\mytitle{$N=\myN$}%
					\ifdefined\Space\empty\else\newcommand{\Space}[1]{#1}\fi%

\ifdefined\mylowerylimit\empty\else\def\mylowerylimit{1e-1}\fi%
\ifdefined\myupperylimit\empty\else\def\myupperylimit{5e2}\fi%
\ifdefined\myN\empty\else\def\myN{1}\fi%
\ifdefined\mydim\empty\else\def\mydim{5000}\fi%
\ifdefined\mypath\empty\else\newcommand{\mypath}{/home/heber/Documents/Experiments/BregmanProjectedSearchDirections}\fi%

\pgfplotsset{
    discard if not/.style 2 args={
        x filter/.code={%
            \edef\tempa{\thisrow{#1}}%
            \edef\tempb{#2}%
            \ifx\tempa\tempb%
            \else%
                \def\pgfmathresult{inf}%
            \fi%
        }
    }
}

\begin{tikzpicture}[]
\begin{loglogaxis}[
	ymin=\mylowerylimit,
	ymax=\myupperylimit,
	grid=major,
	xtick={1.2,1.5,2,3,6,10},
	extra x ticks={1.1},
	extra x tick labels={1.1},
	major tick length={10pt},
	extra x tick style={rotate=0,major tick length=20pt,xtick pos=left, xtick align=center},
	ymajorgrids=true,
	log ticks with fixed point,
	xticklabel shift=2pt,
	legend pos=south east,
	legend columns=1,
	ylabel={runtime [s]},
	xlabel={$\ell_p$-norm}
]
\addplot+[
 	thick,
	discard if not={N}{\myN},
	error bars/.cd,y dir=both,y explicit,
] table[x=p,y=avg_runtime,y error expr={sqrt(\thisrow{var_runtime}-\thisrow{avg_runtime}^2)},col sep=comma]{\mypath /RandomMatrix/data/pnorm-runtimes-Competition_of_search_spaces-dim_\mydim-orthogonal_directions_0.csv};
\addlegendentry{$\Space{U}^{\text{trunc}}_n$}

\addplot+[
 	thick,
	discard if not={N}{\myN},
	error bars/.cd,y dir=both,y explicit,
] table[x=p,y=avg_runtime,y error expr={sqrt(\thisrow{var_runtime}-\thisrow{avg_runtime}^2)},col sep=comma,]{\mypath /RandomMatrix/data/pnorm-runtimes-Competition_of_search_spaces-dim_\mydim-orthogonal_directions_1.csv};
\addlegendentry{$\Space{V}^{\text{trunc}}_n$}

\ifdefined\mytitle{%
\node[rectangle,draw,thick,fill=white] (title) at (1,5) {\mytitle};
}\fi%

\end{loglogaxis}
\end{tikzpicture}
				} \\
				\resizebox{.45\textwidth}{!}{
					\def\myN{2}%
					\edef\mytitle{$N=\myN$}%
					\ifdefined\Space\empty\else\newcommand{\Space}[1]{#1}\fi%

\ifdefined\mylowerylimit\empty\else\def\mylowerylimit{1e-1}\fi%
\ifdefined\myupperylimit\empty\else\def\myupperylimit{5e2}\fi%
\ifdefined\myN\empty\else\def\myN{1}\fi%
\ifdefined\mydim\empty\else\def\mydim{5000}\fi%
\ifdefined\mypath\empty\else\newcommand{\mypath}{/home/heber/Documents/Experiments/BregmanProjectedSearchDirections}\fi%

\pgfplotsset{
    discard if not/.style 2 args={
        x filter/.code={%
            \edef\tempa{\thisrow{#1}}%
            \edef\tempb{#2}%
            \ifx\tempa\tempb%
            \else%
                \def\pgfmathresult{inf}%
            \fi%
        }
    }
}

\begin{tikzpicture}[]
\begin{loglogaxis}[
	ymin=\mylowerylimit,
	ymax=\myupperylimit,
	grid=major,
	xtick={1.2,1.5,2,3,6,10},
	extra x ticks={1.1},
	extra x tick labels={1.1},
	major tick length={10pt},
	extra x tick style={rotate=0,major tick length=20pt,xtick pos=left, xtick align=center},
	ymajorgrids=true,
	log ticks with fixed point,
	xticklabel shift=2pt,
	legend pos=south east,
	legend columns=1,
	ylabel={runtime [s]},
	xlabel={$\ell_p$-norm}
]
\addplot+[
 	thick,
	discard if not={N}{\myN},
	error bars/.cd,y dir=both,y explicit,
] table[x=p,y=avg_runtime,y error expr={sqrt(\thisrow{var_runtime}-\thisrow{avg_runtime}^2)},col sep=comma]{\mypath /RandomMatrix/data/pnorm-runtimes-Competition_of_search_spaces-dim_\mydim-orthogonal_directions_0.csv};
\addlegendentry{$\Space{U}^{\text{trunc}}_n$}

\addplot+[
 	thick,
	discard if not={N}{\myN},
	error bars/.cd,y dir=both,y explicit,
] table[x=p,y=avg_runtime,y error expr={sqrt(\thisrow{var_runtime}-\thisrow{avg_runtime}^2)},col sep=comma,]{\mypath /RandomMatrix/data/pnorm-runtimes-Competition_of_search_spaces-dim_\mydim-orthogonal_directions_1.csv};
\addlegendentry{$\Space{V}^{\text{trunc}}_n$}

\ifdefined\mytitle{%
\node[rectangle,draw,thick,fill=white] (title) at (1,5) {\mytitle};
}\fi%

\end{loglogaxis}
\end{tikzpicture}
				}
				\resizebox{.45\textwidth}{!}{
					\def\myN{6}%
					\edef\mytitle{$N=\myN$}%
					\ifdefined\Space\empty\else\newcommand{\Space}[1]{#1}\fi%

\ifdefined\mylowerylimit\empty\else\def\mylowerylimit{1e-1}\fi%
\ifdefined\myupperylimit\empty\else\def\myupperylimit{5e2}\fi%
\ifdefined\myN\empty\else\def\myN{1}\fi%
\ifdefined\mydim\empty\else\def\mydim{5000}\fi%
\ifdefined\mypath\empty\else\newcommand{\mypath}{/home/heber/Documents/Experiments/BregmanProjectedSearchDirections}\fi%

\pgfplotsset{
    discard if not/.style 2 args={
        x filter/.code={%
            \edef\tempa{\thisrow{#1}}%
            \edef\tempb{#2}%
            \ifx\tempa\tempb%
            \else%
                \def\pgfmathresult{inf}%
            \fi%
        }
    }
}

\begin{tikzpicture}[]
\begin{loglogaxis}[
	ymin=\mylowerylimit,
	ymax=\myupperylimit,
	grid=major,
	xtick={1.2,1.5,2,3,6,10},
	extra x ticks={1.1},
	extra x tick labels={1.1},
	major tick length={10pt},
	extra x tick style={rotate=0,major tick length=20pt,xtick pos=left, xtick align=center},
	ymajorgrids=true,
	log ticks with fixed point,
	xticklabel shift=2pt,
	legend pos=south east,
	legend columns=1,
	ylabel={runtime [s]},
	xlabel={$\ell_p$-norm}
]
\addplot+[
 	thick,
	discard if not={N}{\myN},
	error bars/.cd,y dir=both,y explicit,
] table[x=p,y=avg_runtime,y error expr={sqrt(\thisrow{var_runtime}-\thisrow{avg_runtime}^2)},col sep=comma]{\mypath /RandomMatrix/data/pnorm-runtimes-Competition_of_search_spaces-dim_\mydim-orthogonal_directions_0.csv};
\addlegendentry{$\Space{U}^{\text{trunc}}_n$}

\addplot+[
 	thick,
	discard if not={N}{\myN},
	error bars/.cd,y dir=both,y explicit,
] table[x=p,y=avg_runtime,y error expr={sqrt(\thisrow{var_runtime}-\thisrow{avg_runtime}^2)},col sep=comma,]{\mypath /RandomMatrix/data/pnorm-runtimes-Competition_of_search_spaces-dim_\mydim-orthogonal_directions_1.csv};
\addlegendentry{$\Space{V}^{\text{trunc}}_n$}

\ifdefined\mytitle{%
\node[rectangle,draw,thick,fill=white] (title) at (1,5) {\mytitle};
}\fi%

\end{loglogaxis}
\end{tikzpicture}
				}
				\caption{Runtimes solving with \protect\glostext{SESOP} and two different search spaces for the minimum-norm solution of $\Mat{A}\Vect{x}=\Vect{y}$ with uniformly random $\Mat{A} \in \reel^{1000 \times 5000}$ and $\protect\sym{number-search-directions} \in \{1,2,4,6\}$.}
				\label{fig:solution_iterations-comparingprojected_runtimes-many}
			\end{figure}

			Finally, we look at runtimes for multiple search directions in Figure~\ref{fig:solution_iterations-comparingprojected_runtimes-many}.
			
			We see that for more than one search direction $\sym{number-search-directions}$ the additional cost for the orthogonalization procedure of the search directions at some point starts to outweigh the gain obtained by a reduced number of total iterations. There are two trends behind this: First, \glostext{SESOP} with an unorthogonalized search space also becomes more effective when using multiple search directions, c.\,f.~Table~\ref{tab:solution_iterations-comparingprojected}. Second, the overhead of orthogonalization becomes more costly when $\sym{number-search-directions}$ increases. These two trends work especially in favor of \glostext{SESOP} with $\Space{U}^{\text{trunc}}_n$ for $\sym{norm-X}$ around $2$ where the required iteration steps are comparable for $\sym{number-search-directions} \in \{4,6\}$ and explains why the unorthogonalized procedure is then faster in the total runtimes lacking the additional orthogonalization. This holds for $\sym{norm-X} \in \{1.5,2,3\}$.
			
			Looking at the trend in runtimes for $\Space{V}^{\text{trunc}}_n$ for increasing number of search directions $\sym{number-search-directions}_n$, we see little to no gain. Runtimes for $\sym{norm-X} \leq 2$ generally increase. This trend turns around for towards $\sym{norm-X} > 2$ where the runtimes decrease when using more search directions.
	
		\subsubsection{Connection to \glsentrytext{glos:CG} in Hilbert Space}\label{sec:experiments-toyproblem-connection-CG}
		
			We briefly want to review the $\ell_2$-case, where the Banach space is a Hilbert space. There, equation~\eqref{math:orthogonalized-search-directions-single-direction} is identical to the update used by the Polak-Ribi\`ere~\cite{Polak1969} type of \glostext{CG} method on the normal equation,
			\begin{equation*}
				\Mat{B} \Vect{x} = \Vect{z} \quad \text{with } \Mat{B}=\dual{\Mat{A}}\Mat{A}, \Vect{z} = \dual{\Mat{A}}\Vect{y}.
			\end{equation*} 
			Hence, we expect that the conjugacy property will hold in our case, too. Note that it takes the following form with $\Vect{p}_i = \Mat{A} \Vect{x}_i - \Vect{y}$,
			\begin{equation}
				\skp{ \dual{\Mat{A}} \Vect{p}_i }{ \dual{\Mat{A}} \Vect{p}_j } = 0 \quad \forall i \neq j.
			\end{equation}
			
			Furthermore, we have seen in Table~\ref{tab:solution_iterations-comparingprojected} that there is no difference between a single search direction and multiple search directions as the current gradient direction is made orthogonal to all previous search directions simultaneously.
		
			\begin{figure}[htbp]
				\def\mydim{5000}%
				\def\myseed{426}%
				\def\myorthodir{1}
				\pgfplotsset{height=4cm,width=0.95\textwidth}
				\subfigure[$p=1.2$]{
					\def\mypnorm{1_2}%
					\input{./figures/bregman-distance_iteration}
				} \\
				\subfigure[$p=2$]{
					\def\mypnorm{2}%
					\input{./figures/bregman-distance_iteration}
				} \\
				\subfigure[$p=6$]{
					\def\mypnorm{6}%
					\input{./figures/bregman-distance_iteration}
				}
				\caption{Bregman distance $\protect\sym{distance-bregman}_{\Space{X}} (\Vect{x}_n, \Vect{z})$ over the iteration step $n$ for \protect\glostext{SESOP} with search space $\Space{V}_n^{\text{trunc}}$ for increasing number of search directions $\protect\sym{number-search-directions}$, for various $\ell_{\sym{norm-X}}$ norms with $m=\mydim$, and for a specific random number generator seed.}
				\label{fig:bregman-distance_iteration}
			\end{figure}
	
			In Figure~\ref{fig:bregman-distance_iteration} we give the Bregman distance to the in our case known true solution $\sym{truesolution}$ per iteration step for three different $\ell_{\sym{norm-X}}$ norms, one of them being the $\ell_2$ norm. There, we only look at the orthogonalized search space $\Space{V}_n^{\text{trunc}}$ using metric projections. 
			
			We clearly notice that in the $\ell_2$ case more than a single search direction does not change anything about the minimization. It proceeds (up to numerical precision) in exactly the same manner as if there were only a single search direction spanning the search space. This is precisely what would be expected of a \glostext{CG} method due to the inherent conjugacy property. This is not valid for the other two cases with $p=1.2$ and $p=6$ where more search directions clearly further aid convergence, especially for $p \gg 2$.
			
			Hence, we also see in this numerical example that the subspace methods (with orthogonalization) in fact are an extension of the \glostext{CG} methods to general Banach spaces.
	
	\subsection{Computerized Tomography}\label{sec:experiments-ct}
	
		The computerized tomography problem in 2D aims at reconstruction of the inside of an object from projections. Measurements are for example obtained by passing radiation through a body, whereby their intensity is diminished, proportional to the passed length and density of the body $f(x): [0,1]^2 \rightarrow \reel^+_0$. This decrease is measured over a $a$ angles and $s$ shifts of radiation source and detectors. This measurement matrix is usually called \emph{sinogram}. We follow \cite[\citesection~7.7]{Hansen2010a} in the brief introduction to the discretized setting.

		The measurement rays are parametrized as
		\begin{equation}\label{math:ray-parametrization}
			t^i(\tau) = t^{i,0} + \tau d^i,
		\end{equation}
		where $d$ is the directional vector of the $i$th ray with $i=1, \ldots, s\cdot a$.
		
		Then, the dampening of the $i$th ray can be written as the \emph{Radon transform},
		\begin{equation}\label{math:radon-transform}
			b_i = \int^{\infty}_{-\infty} f \left (t^i(\tau) \right ) d\tau, \quad i=1, \ldots, s\cdot a,
		\end{equation}
		i.\,e.~we integrate $d\tau$ along the line $t^i(\tau)$.
		
		The problem can be discretized with a pixel basis,
		\begin{equation}\label{math:pixel-basis}
			\chi_{kl} (x) = \begin{cases} 1, &x \in [h \cdot (k-1), h \cdot k] \times [h \cdot (l-1), h \cdot l] \\ 0, &\text{else} \end{cases}, \quad \forall k,l=1,\ldots,n,
		\end{equation}
		where we constrain the absorption coefficient $f(x) = \sum_{kl} f_{kl} \chi_{kl} (x)$ to be piecewise constant with $h=\tfrac 1 n$.
		
		For the discretization of \eqref{math:radon-transform} we simply need to count the length $\Delta L^{(i)}_{kl}$ of each ray $t^i$ in each pixel $\chi_{kl}$ of the basis~\eqref{math:pixel-basis} and obtain
		\begin{equation}\label{math:radon-transform-discretized}
			b_i = \sum^n_{k,l=1} f_{kl} \Delta L^{(i)}_{kl} \quad \text{for } i=1,\ldots,s\cdot a.
		\end{equation}
		
		If we vectorize the matrix object $f_{kl}$ to become the vector $x_j$ with $j=(l-1)n+k$, then we obtain
		\begin{equation}\label{math:radon-transform-discretized-vectorized}
			b_i = \sum_{j} a_{ij} x_j, \quad i=1,\ldots,s\cdot a.
		\end{equation}
		
		Note that the matrix $A$ is sparse what we use in the implementation.

			\begin{figure*}[htbp]
				\def\myprefix{SESOP}%
				\def\myangles{60}%
				\def\myoffsets{61}%
				\def\myp{2}%
				\def\mypixels{41}%
				\begin{center}
				\subfigure[Phantom\label{fig:phantom}]{
					\includegraphics[width=0.45\textwidth]{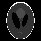} 
				} \goodgap
				\subfigure[Sinogram\label{fig:sinogram}]{
					\includegraphics[width=0.458\textwidth]{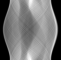} 
				}
				\end{center}
				\caption{Shepp-Logan phantom and its analytically obtained sinogram.}
				\label{fig:problem}
			\end{figure*}	

			We look at the standard Shepp-Logan phantom, see Figure~\ref{fig:phantom}, where we obtain the measurements, see Figure~\ref{fig:sinogram}, by using the known analytical Radon transform of ellipses. We look at only a single exemplary discretization with number of pixels $n=41$, number of shifts $s=61$, and number of discrete angles $a=60$. We deliberately choose a coarser image resolution together with a higher number of measurements to allow for a high-quality reconstruction and clearly discernable artifacts if there are any. 
			
			Note that we additionally project the solution onto the range of the matrix $A$ such that the true solution can be computed by the method. We return to this point in Section \ref{sec:experiments-ct-with-noise}.

			Here, we again want to compare \glostext{SESOP} using the truncated search space $\Space{U}^{\text{trunc}}_n$ and also the orthogonalized search space $\Space{V}^{\text{trunc}}_n$. They have the following dimensions: $\abs{\Space{U}^{\text{trunc}}_n} = 2$ and $\abs{\Space{V}^{\text{trunc}}_n}=1$. We stop the methods either after $500$ iterations or when the absolute residual has decreased below $10^{-2}$. We use $\ell_{\sym{norm-X}}$-spaces with $\sym{norm-X} \in \{1.1,1.2,1.5,2\}$ for $\Space{X}$ and $\ell_2$ for $\Space{Y}$. For the line-search problem we use at most $20$ iterations. The power type of the duality mappings is always set to $2$.
		
		\subsubsection{Exact data}\label{sec:experiments-ct-without-noise}
			First of all, we look at the problem without noise and only in the case of $\Space{X} = \ell_2$, i.\,e.~the Hilbert space setting.
		
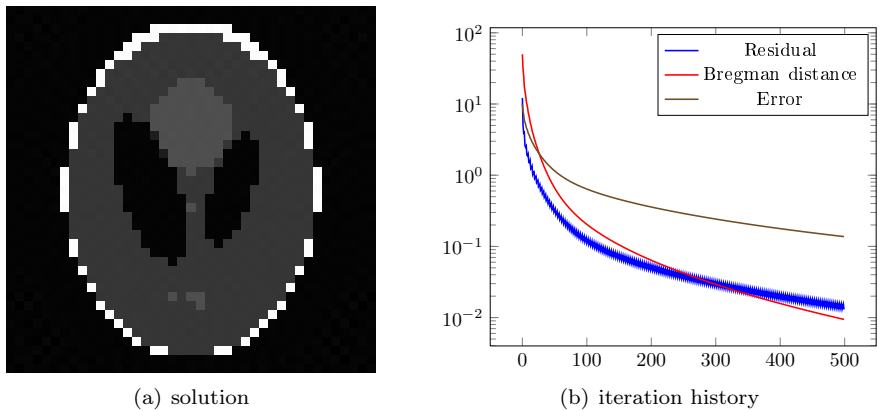
\begin{figure*}[htbp]
				\def\myprefix{SESOP}%
				\def\myangles{60}%
				\def\myoffsets{61}%
				\def\myp{2}%
				\def\mypixels{41}%
				\begin{center}
				\subfigure[solution\label{fig:solution-SESOP}]{
					\includegraphics[width=0.4\textwidth]{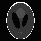} 
				} \goodgap
				\subfigure[iteration history\label{fig:residual-SESOP}]{
					\resizebox{0.5\textwidth}{!}{
						\ifdefined\myprefix\empty\else\def\myprefix{SESOP}\fi%
\ifdefined\myangles\empty\else\def\myangles{40}\fi%
\ifdefined\myoffsets\empty\else\def\myoffsets{41}\fi%
\ifdefined\myp\empty\else\def\myp{2}\fi%
\ifdefined\mypixels\empty\else\def\mypixels{41}\fi%
\ifdefined\mypath\empty\else\newcommand{\mypath}{/home/heber/Documents/Experiments/BregmanProjectedSearchDirections}\fi%

\begin{tikzpicture}[]
\begin{semilogyaxis}
\addplot+[mark=none,thick] table[x=iteration,y=residual, col sep=comma]{\mypath /ComputerizedTomography/data/residual-\myprefix-num_angles_\myangles-num_offsets_\myoffsets-num_pixels_\mypixels-p_\myp.dat};
\addlegendentry{Residual}
\addplot+[mark=none,thick] table[x=iteration,y=bregman_distance, col sep=comma]{\mypath /ComputerizedTomography/data/residual-\myprefix-num_angles_\myangles-num_offsets_\myoffsets-num_pixels_\mypixels-p_\myp.dat};
\addlegendentry{Bregman distance}
\addplot+[mark=none,thick] table[x=iteration,y=error, col sep=comma]{\mypath /ComputerizedTomography/data/residual-\myprefix-num_angles_\myangles-num_offsets_\myoffsets-num_pixels_\mypixels-p_\myp.dat};
\addlegendentry{Error}
\end{semilogyaxis}
\end{tikzpicture}
					}
				}
				\end{center}
				\caption{Recovered solution and residual, error, and Bregman distance histories using SESOP with $\abs{\Space{U}^{\text{trunc}}_n} = 2$ directions.}
				\label{fig:problem-SESOP}
			\end{figure*}

			In Figure~\ref{fig:problem-SESOP} we give the recovered solution and the error, residual, and Bregman distance histories for \glostext{SESOP} with the truncated search space $\Space{U}^{\text{trunc}}_n$ and two search directions.

			We see that \glostext{SESOP} converges but stops at $n=500$ iterations where the absolute residual is still slightly larger than $10^{-2}$. The overall runtime is $2.4$~seconds. Note that the convergence is monotone only in the Bregman distance and not in the residual, c.\,f.~Thm~\ref{thm:convergence-truncated-orthogonalized-searchspace}.

\begin{figure*}[htbp]
				\def\myprefix{metricSESOP}%
				\def\myangles{60}%
				\def\myoffsets{61}%
				\def\myp{2}%
				\def\mypixels{41}%
				\begin{center}
				\subfigure[solution\label{fig:solution-metricSESOP}]{
					\includegraphics[width=0.4\textwidth]{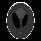} 
				} \goodgap
				\subfigure[iteration history\label{fig:residual-metricSESOP}]{
					\def\myprefix{metricSESOP}%
					\resizebox{0.5\textwidth}{!}{
						\ifdefined\myprefix\empty\else\def\myprefix{SESOP}\fi%
\ifdefined\myangles\empty\else\def\myangles{40}\fi%
\ifdefined\myoffsets\empty\else\def\myoffsets{41}\fi%
\ifdefined\myp\empty\else\def\myp{2}\fi%
\ifdefined\mypixels\empty\else\def\mypixels{41}\fi%
\ifdefined\mypath\empty\else\newcommand{\mypath}{/home/heber/Documents/Experiments/BregmanProjectedSearchDirections}\fi%

\begin{tikzpicture}[]
\begin{semilogyaxis}
\addplot+[mark=none,thick] table[x=iteration,y=residual, col sep=comma]{\mypath /ComputerizedTomography/data/residual-\myprefix-num_angles_\myangles-num_offsets_\myoffsets-num_pixels_\mypixels-p_\myp.dat};
\addlegendentry{Residual}
\addplot+[mark=none,thick] table[x=iteration,y=bregman_distance, col sep=comma]{\mypath /ComputerizedTomography/data/residual-\myprefix-num_angles_\myangles-num_offsets_\myoffsets-num_pixels_\mypixels-p_\myp.dat};
\addlegendentry{Bregman distance}
\addplot+[mark=none,thick] table[x=iteration,y=error, col sep=comma]{\mypath /ComputerizedTomography/data/residual-\myprefix-num_angles_\myangles-num_offsets_\myoffsets-num_pixels_\mypixels-p_\myp.dat};
\addlegendentry{Error}
\end{semilogyaxis}
\end{tikzpicture}
					}
				}
				\end{center}
				\caption{Recovered solution and residual, error, and Bregman distance histories using orthogonalized SESOP with $\abs{\Space{V}^{\text{trunc}}_n}=1$ direction.}
				\label{fig:problem-metricSESOP}
			\end{figure*}
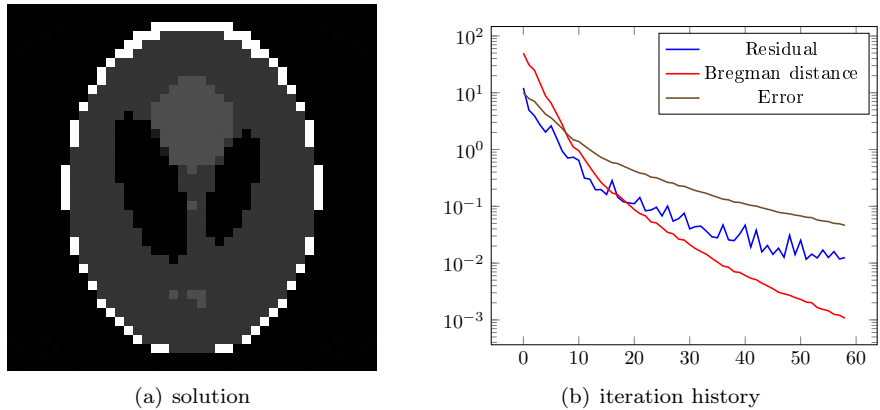

	
			Next, we look at SESOP using the orthogonalized search space $\Space{V}^{\text{trunc}}_n$ with a single search direction\footnote{As mentioned before, in an $\ell_2$ space more search directions do not improve performance when using an orthogonalized search space.}. We obtain the results as given in figure~\ref{fig:problem-metricSESOP}.
			
			Here, we notice that the iteration stops at roughly $n=60$ when reaching an absolute residual of $10^{-2}$. Naturally, this translates to a faster runtime compared to the unorthogonalized search space of only $0.36$~seconds.

			In the former case the reconstructed image retains some very slight artefacts, in the latter case no artefacts are visible.			

		\subsubsection{Noisy data}\label{sec:experiments-ct-with-noise}
		
			It is called \emph{inverse crime}, see \cite[\citesection~7.2]{Hansen2010a}, if the same discretization has been used for both the operator and the right-hand side. Results may look suspiciously good in this case, i.\,e.~no mismatch between (real) data and model is revealed. We definitely committed this crime by projecting our measurements onto the range of the matrix $A$. Now, we additionally disturb the right-hand side $\Vect{y}$, obtained from projecting the Shepp-Logan phantom $\sym{truesolution}$ with the matrix $\Vect{y} = \Mat{A} \sym{truesolution}$, with noise of a known level $\delta$ to become $\Vect{y}^{\delta}$. This ensures that $\norm{ \Vect{y} - \Vect{y}^{\delta} }$ is known to us.
						
			To this end, we construct a random vector $\Vect{n} \in [-1,1]^{\sym{dimension-solution}}$ and set $\widetilde{\Vect{y}} = \Vect{y} + \delta \frac{\norm{\Vect{y}}} {\norm{\Vect{n}}} \Vect{n}$ with the projected right-hand side $\Vect{y}=\Mat{A}\sym{truesolution}$ using the true solution $\sym{truesolution}$, which is the Shepp-Logan phantom. We use a noise level of $\delta=0.01$. We stop the iteration via a relative residual threshold of $0.01$ with a discrepancy parameter of $3$, i.\,e.~the relative residual threshold becomes $0.03$. 
	
			As a perfect reconstruction is in the noisy case no longer possible, we increase the dimensions of the measurement matrix to number of offsets of $s=81$ and number of angles to $a=80$. Also we increase the number of pixels $p=81$ of the reconstructed image. We still remain in the Hilbert space setting $\Space{X} = \ell_2$ for the moment. In the next section we look at other $\ell_{\sym{norm-X}}$ spaces.
	
			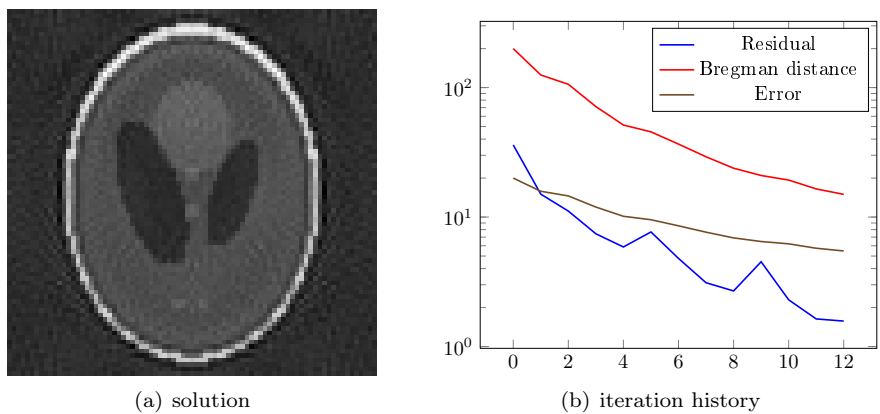
\begin{figure*}[htbp]
				\def\myprefix{inverseCrime}%
				\def\myangles{80}%
				\def\myoffsets{81}%
				\def\myp{2}%
				\def\mypixels{81}%
				\def\myorthodir{1}%
				\begin{center}
				\subfigure[solution\label{fig:solution-inversecrime}]{
					\includegraphics[width=0.4\textwidth]{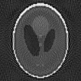} 
				} \goodgap
				\subfigure[iteration history\label{fig:residual-inversecrime}]{
					\resizebox{0.5\textwidth}{!}{
						\ifdefined\myangles\empty\else\def\myangles{80}\fi%
\ifdefined\myoffsets\empty\else\def\myoffsets{81}\fi%
\ifdefined\myp\empty\else\def\myp{2}\fi%
\ifdefined\mypixels\empty\else\def\mypixels{81}\fi%
\ifdefined\myorthodir\empty\else\def\myorthodir{1}\fi%
\ifdefined\mypath\empty\else\newcommand{\mypath}{/home/heber/Documents/Experiments/BregmanProjectedSearchDirections}\fi%

\begin{tikzpicture}[]
\begin{semilogyaxis}
\addplot+[mark=none,thick] table[x=iteration,y=residual, col sep=comma]{\mypath /ComputerizedTomography/data/residual-inverseCrime-num_angles_\myangles-num_offsets_\myoffsets-num_pixels_\mypixels-orthogonal_directions_\myorthodir-p_\myp.dat};
\addlegendentry{Residual}
\addplot+[mark=none,thick] table[x=iteration,y=bregman_distance, col sep=comma]{\mypath /ComputerizedTomography/data/residual-inverseCrime-num_angles_\myangles-num_offsets_\myoffsets-num_pixels_\mypixels-orthogonal_directions_\myorthodir-p_\myp.dat};
\addlegendentry{Bregman distance}
\addplot+[mark=none,thick] table[x=iteration,y=error, col sep=comma]{\mypath /ComputerizedTomography/data/residual-inverseCrime-num_angles_\myangles-num_offsets_\myoffsets-num_pixels_\mypixels-orthogonal_directions_\myorthodir-p_\myp.dat};
\addlegendentry{Error}
\end{semilogyaxis}
\end{tikzpicture}
					}
				}
				\end{center}
				\caption{Recovered solution and residual, error, and Bregman distance histories using orthogonalized SESOP with $\abs{\Space{V}^{\text{trunc}}_n}=1$ direction in the presence of noise of level $0.01$ and using a stopping criterion with $\tau=3$.}
				\label{fig:problem-inversecrime}
			\end{figure*}

			In Figure~\ref{fig:problem-inversecrime} both the reconstructed image and iteration history with residual, Bregman distance, and error with respect to the true solution is shown. We obtain good results with respect to the noise level employed. The runtime is $0.34$~seconds.
			
			We conclude that we do not commit any inverse crime and that the proposed method is indeed working and implemented sensibly.
			
		
		\subsubsection{Small $\ell_p$ Norms}\label{sec:experiments-ct-procedure-other-norms}
		
			Finally, we also look at other norms than $\ell_2$, namely small $\ell_p$ norms with $\sym{norm-X} \in \{1.1, 1.2, 1.5\}$. Here, to allow for a high-quality reconstruction we again use $s=61$, $a=60$ and $p=41$.

			\begin{figure}[htbp]
				\def\myangles{60}%
				\def\myoffsets{61}%
				\def\mypixels{41}%
				\pgfplotsset{height=6cm,width=0.95\textwidth}
				\begin{center}
				\input{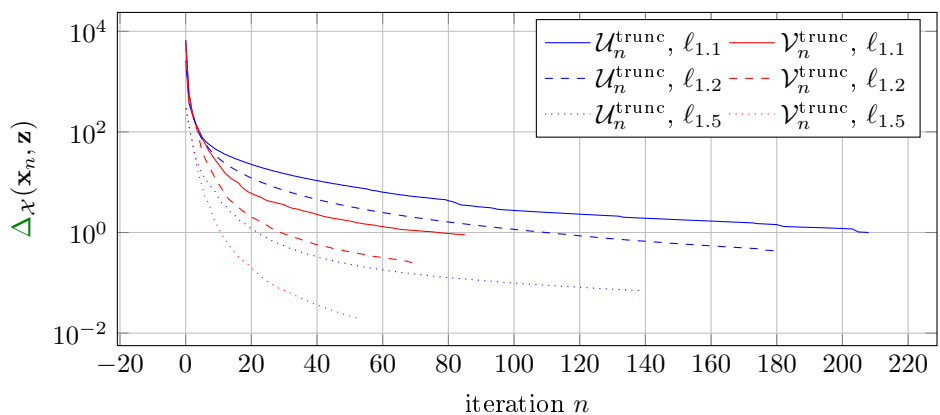}
				\end{center}
				\caption{Bregman distance $\protect\sym{distance-bregman}_{\Space{X}} (\Vect{x}_n, \Vect{z})$ over the iteration step $n$ for \protect\glostext{SESOP} with unorthogonalized search space $\Space{U}_n^{\text{trunc}}$ and orthogonalized search space $\Space{V}_n^{\text{trunc}}$ in comparison on the norms $\ell_{1.1}$, $\ell_{1.2}$, and $\ell_{1.5}$.}
				\label{fig:bregman-distance_iteration-smallerlp}
			\end{figure}

			In Figure~\ref{fig:bregman-distance_iteration-smallerlp} we compare the decrease in Bregman distance for unorthogonalized and orthogonalized search spaces for various small $\ell_p$ norms. We see that over all tested small $\ell_p$ norms the orthogonalized search space $\Space{V}_n^{\text{trunc}}$ requires about three times fewer iterations to reach the stopping criterion. This makes it also faster in the overall runtimes, e.\,g.~for $\sym{norm-X}=1.1$ we obtain $3.15$~seconds without orthogonalization and $2.6$~seconds when employing $\Space{V}_n^{\text{trunc}}$.
			Note that because of different $\ell_{\sym{norm-X}}$ used to measure the residual stopping criterion, the method also terminates at different Bregman distances.
			
			\begin{figure}[htbp]
				\def\myprefix{inverseCrime}%
				\def\myangles{80}%
				\def\myoffsets{81}%
				\def\mypixels{127}%
				\def\myorthodir{1}%
				\begin{center}
				\subfigure[$\ell_{1.1}$\label{fig:reconstructions-smallerlp-1_1}]{
					\def\myp{1_1}%
					\includegraphics[width=0.2\textwidth]{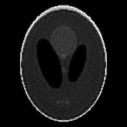} 
				} \goodgap
				\subfigure[$\ell_{1.2}$\label{fig:reconstructions-smallerlp-1_2}]{
					\def\myp{1_2}%
					\includegraphics[width=0.2\textwidth]{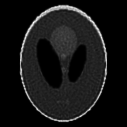} 
				} \goodgap
				\subfigure[$\ell_{1.5}$\label{fig:reconstructions-smallerlp-1_5}]{
					\def\myp{1_5}%
					\includegraphics[width=0.2\textwidth]{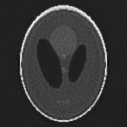} 
				} \goodgap
				\subfigure[$\ell_{2}$\label{fig:reconstructions-smallerlp-2}]{
					\def\myp{2}%
					\includegraphics[width=0.2\textwidth]{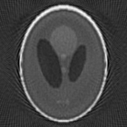} 
				}
				\end{center}
				\caption{Reconstructed Shepp-Logan phantom in the presence of noise of level $\delta=0.01$ and using small $\ell_{\sym{norm-X}}$ norms with $\sym{norm-X} \in \{1.1,1.2,1.5,2\}$.}
				\label{fig:reconstructions-smallerlp}
			\end{figure}

			Last of all, we show in Figure~\ref{fig:reconstructions-smallerlp} larger reconstructed images for $p=127$ pixels with a noise of $\delta=0.01$ while using same number of angles $a$ and shifts $s$ than before. This is to elucidate the effect of the different $\ell_p$ norms in the presence of noise.

			We observe that for smaller $\ell_{\sym{norm-X}}$ both contrast and also noise of the image is enhanced. Especially, artifacts surrounding the reconstructed phantom in the $\ell_2$ case are absent for $\ell_{1.1}$. On the other hand, noisy speckles are more present in the reconstructed image using the $\ell_{1.1}$ norm. 

			
\section{Conclusions}\label{sec:conclusions}

	Based on the previous work of \cite{Schopfer2007}, we have proposed a semi-orthogonalized set of search directions in Banach spaces using metric projections. Using search spaces consisting only of the last Landweber descent directions modified by this orthogonalization, we have shown that the \glosfirst{SESOP} method converges weakly and that for the Hilbert space case the procedure coincides with the Polak-Ribi\`ere \glosfirst{CG} method applied to the normal equation. Hence, the subspace methods are a natural extension to general Banach spaces of \glostext{CG} methods, that are known to work very efficiently in Hilbert spaces. 
	
	Numerical experiments have shown fast convergence for both, an inverse toy problem consisting of a uniformly distributed random matrix and right-hand side on various $\ell_{\sym{norm-X}}$ spaces as well as for the inverse problem of 2D computerized tomography using the discretized backprojection as operator and a Radon transformed Shepp-Logan phantom as right-hand side. In every case the orthogonalized truncated search space clearly outperforms the truncated search space used in \cite{Schopfer2008}, both in required iterations and runtime.

	Note that although we have investigated the case of a noisy right-hand side in the experiments, convergence for this case has not been done yet.

	As an outlook, we would like to remind that there is a whole zoo of \glostext{CG} variants, see \cite[\citepage~98]{Andrei2009}. It would be very insightful to find more connections between a specific variant and a choice of (semi-orthogonalized) search spaces for the \glostext{SESOP} method. Proving the regularization property of the method is another object of future research.

\section*{Acknowledgements}\label{sec:acknowledgements}

The authors acknowledge funding by the German ministry of Education and Research (BMBF) under the project HYPERMATH (05M2013).


\bibliographystyle{siamplain}	
 \bibliography{small_library}

\end{document}